\documentclass[a4paper]{amsart}

\usepackage{amssymb}

\makeatletter
\def\section{\@startsection{section}{1}%
  \z@{.7\linespacing\@plus\linespacing}{.5\linespacing}%
  {\normalfont\bfseries\centering}}
\def\@secnumfont{\bfseries}
\makeatother

\def\frak{\mathfrak}
\def\Bbb{\mathbb}
\def\Cal{\mathcal}
\def\sideremark#1{\ifvmode\leavevmode\fi\vadjust{\vbox to0pt{\vss%
  \hbox to 0pt{\hskip\hsize\hskip1em%
  \vbox{\hsize3cm\tiny\raggedright\pretolerance10000%
  \noindent #1\hfill}\hss}\vbox to8pt{\vfil}\vss}}}%
\def\leftsideremark#1{\ifvmode\leavevmode\fi
\vadjust{\vbox to0pt{\vss%
 \hbox to 0pt{\hskip-4.5cm%
 \vbox{\hsize4cm\tiny\raggedright\pretolerance10000%
 \noindent \hfill#1}\hss}\vss}}}%

\swapnumbers

\newcommand{\al}{\alpha}

\newcommand{\ka}{\kappa}

\newcommand{\om}{\omega}
\newcommand{\ph}{\varphi}
\newcommand{\ps}{\psi}

\newcommand{\ze}{\zeta}

\newcommand{\La}{\Lambda}
\newcommand{\Ph}{\Phi}
\newcommand{\Ps}{\Psi}

\newcommand{\Om}{\Omega}

\newcommand{\ad}{\operatorname{ad}}
\newcommand{\Ad}{\operatorname{Ad}}
\newcommand{\Adb}{\operatorname{\underline{Ad}}}

\renewcommand{\exp}{\operatorname{exp}}
\newcommand{\id}{\operatorname{id}}

\newcommand{\im}{\operatorname{im}}
\newcommand{\tr}{\operatorname{tr}}

\newcommand{\sgn}{\operatorname{sgn}}

\let\[=\llbracket
\let\]=\rrbracket
\let\x=\times

\def\g{\frak g}
\newcommand{\fg}{{\frak g}}

\def\({\big(}
\def\){\big)}
\def\R{\Bbb R}

\def\G{{\Cal G}}
\newcommand{\tcg}{{\tilde{\Cal G}}}
\newcommand{\tg}{{\tilde{\frak g}}}
\newcommand{\tp}{{\tilde{\frak p}}}

\def\.{\hbox to5pt{\hss$\cdot$\hss}}

\newtheorem*{prop*}{Proposition}

\newtheorem*{thm*}{Theorem}

\newtheorem*{lemma*}{Lemma}

\newtheorem*{cor*}{Corollary}

\renewcommand{\o}{\circ}

\begin{document}

\title{On the Geometry of Chains}
\author{Andreas \v Cap\\ Vojt\v ech \v Z\'adn\'\i k}
\thanks{First author supported by project
P15747--N05 of the Fonds zur F\"orderung der wissenschaftlichen
Forschung (FWF). Second author supported at different times by the Junior 
Fellows program of the Erwin Schr\"odinger Institute (ESI) and by the grant 
201/05/2117 of the Czech Science Foundation (GA\v CR).
Discussions with Boris Doubrov have been very helpful.} 

\address{A.\v C: Fakult\"at f\"ur Mathematik, Universit\"at Wien, 
Nordbergstra\ss e 15, A--1090 Wien, Austria and International 
Erwin Schr\"odinger Institute for Mathematical Physics, Boltzmanngasse
9, A--1090 Wien, Austria\newline V.\v Z: International 
Erwin Schr\"odinger Institute for Mathematical Physics, Boltzmanngasse
9, A--1090 Wien, Austria and Faculty of Education, Masaryk University, 
Po\v r\'\i\v c\'\i{} 31, 60300 Brno, Czech Republic}
\email{}
\begin{abstract}
The chains studied in this paper generalize Chern--Moser chains for CR
structures. They form a distinguished family of one dimensional
submanifolds in manifolds endowed with a parabolic contact
structure. Both the parabolic contact structure and the system of
chains can be equivalently encoded as Cartan geometries (of different
types). The aim of this paper is to study the relation between these
two Cartan geometries for Lagrangean contact structures and partially
integrable almost CR structures.

We develop a general method for extending Cartan geometries which
generalizes the Cartan geometry interpretation of Fefferman's
construction of a conformal structure associated to a CR
structure. For the two structures in question, we show that the Cartan
geometry associated to the family of chains can be obtained in that
way if and only if the original parabolic contact structure is torsion
free. In particular, the procedure works exactly on the subclass of
(integrable) CR structures.

This tight relation between the two Cartan geometries leads to an
explicit description of the Cartan curvature associated to the family
of chains. On the one hand, this shows that the homogeneous models for
the two parabolic contact structures give rise to examples of
non--flat path geometries with large automorphism groups. On the other
hand, we show that one may (almost) reconstruct the underlying torsion
free parabolic contact structure from the Cartan curvature associated
to the chains. In particular, this leads to a very conceptual proof of
the fact that chain preserving contact diffeomorphisms are either
isomorphisms or anti--isomorphisms of parabolic contact structures.
\end{abstract}
\date{May 26, 2005}
\maketitle

\section{Introduction}
Parabolic contact structures are a class of geometric structures
having an underlying contact structure. They admit a canonical normal
Cartan connection corresponding to a contact grading of a simple Lie
algebra. The best known examples of such structures are
non--degenerate partially integrable almost CR structures of
hypersurface type. The construction of the canonical Cartan connection
is due to Chern and Moser (\cite{C-M}) for the subclass
of CR structures, and to Tanaka  (\cite{Tan76}) in general. 

In the approach of Chern and Moser, a central role is played by a
canonical class of unparametrized curves called \textit{chains}. For
each point $x$ and each direction $\xi$ at $x$, which is transverse to
the contact distribution, there is a unique chain through $x$ in direction
$\xi$. In addition, each chain comes with a projective class of
distinguished parametrizations. The notion of chains easily
generalizes to arbitrary parabolic contact structures, and the chains
are easy to describe in terms of the Cartan connection.

A \textit{path geometry} on a smooth manifold $M$ is given by a smooth
family of unparametrized curves on $M$ such that for each $x\in M$ and
each direction $\xi$ at $x$ there is a unique curve through $x$ in
direction $\xi$. The best way to encode this structure is to pass to
the projectivized tangent bundle $\Cal PTM$, the space of all lines in
$TM$. Then a path geometry is given by a line subbundle in the tangent
bundle of $\Cal PTM$ with certain properties, see \cite{Douglas} and
\cite{Grossman} for a modern presentation. It turns out that these
structures are equivalent to regular normal Cartan geometries of a
certain type, which fall under the general concept of parabolic
geometries, see section 4.7 of \cite{C-tw}.

In the description as a Cartan geometry, path geometries immediately
generalize to open subsets of the projectivized tangent bundle. In
particular, given a manifold $M$ endowed with a parabolic contact
structure, the chains give rise to a path geometry on the open subset
$\Cal P_0TM\subset\Cal PTM$ formed by all lines transversal to the
contact subbundle. The general question addressed in this paper is how
to describe the resulting Cartan geometry on $\Cal P_0TM$ in terms of
the original Cartan geometry on $M$.  We study this in detail in the
case of Lagrangean contact structures and, in the end, briefly indicate
how to deal with partially integrable almost CR structures, which can
be viewed as a different real form of the same complex geometric
structure.

The first observation is that $\Cal P_0TM$ can be obtained as a
quotient of the Cartan bundle $\Cal G\to M$ obtained from the
parabolic contact structure. More precisely, there is a subgroup
$Q\subset P$ such that $\Cal P_0TM\cong\Cal G/Q$. In particular, $\Cal
G$ is a principal $Q$--bundle over $\Cal P_0TM$ and the canonical
Cartan connection $\om\in\Om^1(\Cal G,\frak g)$ associated to the
parabolic contact structure can be also viewed as a Cartan connection
on $\Cal G\to\Cal P_0TM$. The question then is whether the canonical
Cartan geometry $(\tcg\to\Cal P_0TM,\tilde\om)$ determined by the path
geometry of chains can be constructed directly from $(\Cal G\to\Cal
P_0TM,\om)$.

To attack this problem, we study a class of extension functors mapping
Cartan geometries of some type $(G,Q)$ to Cartan geometries of another
type $(\tilde G,\tilde P)$. These functors have the property that
there is a homomorphism between the two Cartan bundles, which relates
the two Cartan connections. We show that in order to obtain such a functor, one
needs a homomorphism $i:Q\to\tilde P$ (which we assume to be
infinitesimally injective) and a linear map $\al:\frak g\to\tg$ which
satisfy certain compatibility conditions. There is a simple notion of
equivalence for such pairs and equivalent pairs lead to naturally
isomorphic extension functors. 

There is a particular simple source of pairs $(i,\al)$ leading to
extension functors as above. Namely, one may start from a homomorphism
$G\to\tilde G$ and take $i$ the restriction to $Q$ and $\al$ the
induced homomorphism of Lie algebras. In a special case, this leads to
the Cartan geometry interpretation of Fefferman's construction of a
canonical conformal structure on a circle bundle over a CR
manifold. 

One can completely describe the effect of the extension functor
associated to a pair $(i,\al)$ on the curvature of the Cartan
geometries. Apart from the curvature of the original geometry, also
the deviation from $\al$ being a homomorphism of Lie algebras enters
into the curvature of the extended Cartan geometry.

An important feature of the special choice for $(\tilde G,\tilde P)$
that we are concerned with, is a uniqueness result for such extension
functors. We show (see Theorem \ref{3.4}) that if the extension
functor associated to a pair $(i,\al)$ maps locally flat geometries of
type $(G,Q)$ to regular normal geometries of type $(\tilde G,\tilde
P)$, then the pair $(i,\al)$ is already determined uniquely up to
equivalence. For the two parabolic contact structures studied in this
paper, we show that there exist appropriate pairs $(i,\al)$ in
\ref{3.5} and \ref{5.2}.

In both cases, the resulting extension functor does \textit{not}
produce the canonical Cartan geometry associated to the path geometry
of chains in general. We show that the canonical Cartan connection is
obtained if and only if the original parabolic contact geometry is
torsion free. For a Lagrangean contact structure this means that the
two Lagrangean subbundles are integrable, while it is the usual
integrability condition for CR structures. This ties in nicely with
the Fefferman construction, where one obtains a conformal
structure for arbitrary partially integrable almost CR structures, but
the normal Cartan connection is obtained by equivariant extension if
and only if the structure is integrable (and hence CR).

Finally, we discuss applications of our construction, which are based
on an ana\-lysis of the curvature of the canonical Cartan connection
associated to the path geometry of chains. We show that chains never
are geodesics of a connection, and they give rise to a torsion free
path geometry if and only if the original parabolic contact structure
is locally flat. Then we show that the underlying parabolic contact
structure can be almost reconstructed from the harmonic curvature of
the path geometry of chains. In particular, this leads to a very
conceptual proof of the fact that a contact diffeomorphism which maps
chains to chains must (essentially) preserve the original torsion free
parabolic contact structure.

\section{Parabolic contact structures, chains, and path geometries}
\label{2}
In this section, we will discuss the concepts of chains and the
associated path geometry for a parabolic contact structure, focusing
on the example of Lagrangean contact structures. We only briefly
indicate the changes needed to deal with general parabolic contact
structures. 

\subsection{Lagrangean contact structures}\label{2.1}
The starting point to define a parabolic contact structure is a simple
Lie algebra $\frak g$ endowed with a \textit{contact grading}, i.e.~a
vector space decomposition $\frak g=\frak g_{-2}\oplus\frak
g_{-1}\oplus\frak g_0\oplus\frak g_1\oplus\frak g_2$ such that $[\frak
g_i,\frak g_j]\subset\frak g_{i+j}$, $\frak g_{-2}$ has real dimension
one, and the bracket $\frak g_{-1}\x\frak g_{-1}\to\frak g_{-2}$ is
non--degenerate. It is known that such a grading is unique up to an
inner automorphism and it exists for each non--compact non--complex
real simple Lie algebra except $\frak{sl}(n,\Bbb H)$,
$\frak{so}(n,1)$, $\frak{sp}(p,q)$, one real form of $E_6$ and one of
$E_7$, see section 4.2 of \cite{Yam}.

Here we will mainly be concerned with the contact grading of
$\g=\frak{sl}(n+2,\R)$, corresponding to the following block
decomposition with blocks of size $1$, $n$, and $1$:
$$
\begin{pmatrix}
\g_0&\g^L_1&\g_2\\\g^L_{-1}&\g_0&\g^R_1\\\g_{-2}&\g^R_{-1}&\g_0
\end{pmatrix}.
$$
We have indicated the splittings $\frak g_{-1}=\frak
g_{-1}^L\oplus\frak g_{-1}^R$ respectively $\frak g_1=\g_1^L\oplus
g_1^R$, which are immediately seen to be $\fg_0$--invariant. Further, 
the subspaces $\g_{-1}^L$ and $\g_{-1}^R$ of $\g_{-1}$ are isotropic
for $[\ ,\ ]:\g_{-1}\x\g_{-1}\to\g_{-2}$. 

Put $G:=PGL(n+2,\Bbb R)$, the quotient of $GL(n+2,\Bbb R)$ by its
center. We will view $G$ as the quotient of the group of matrices
whose determinant has modulus one by the two element subgroup
generated by $\pm\id$ and work with representative matrices. The group
$G$ always has Lie algebra $\fg$. For odd $n$, one can identify $G$
with $SL(n+2,\Bbb R)$. For even $n$, $G$ has two connected components,
 and the component containing the identity is $PSL(n+2,\Bbb R)$. 
 
 By $G_0\subset P\subset G$ we denote the subgroups formed by matrices
 which are block diagonal respectively block upper triangular with
 block sizes $1$, $n$, and $1$. Then the Lie algebras of $G_0$ and $P$
 are $\fg_0$ respectively $\frak p:=\g_0\oplus\g_1\oplus\g_2$. For
 $g\in G_0$, the map $\Ad(g):\fg\to\fg$ preserves the grading while
 for $g\in P$ one obtains $\Ad(g)(\g_i)\in\g_i\oplus\dots\oplus\g_2$
 for $i=-1,\dots,2$. This can be used as an alternative
 characterization of the two subgroups. The reason for the choice of
 the specific group $G$ with Lie algebra $\fg$ is that the adjoint
 action identifies $G_0$ with the group of all automorphisms of the
 graded Lie algebra $\fg_{-2}\oplus\fg_{-1}$ which in addition
 preserve the decomposition $\fg_{-1}=\fg_{-1}^L\oplus\fg_{-1}^R$.

Let $M$ be a smooth manifold of dimension $2n+1$ and let $H\subset TM$
be a subbundle of corank one. The Lie bracket of vector fields
induces a tensorial map $\Cal L:\La^2H\to TM/H$, and that $H$ is
called a \textit{contact structure} on $M$ if this map is
non--degenerate. A \textit{Lagrangean contact structure} on $M$ is a
contact structure $H\subset TM$ together with a fixed decomposition
$H=L\oplus R$ such that each of the subbundles is isotropic with
respect to $\Cal L$. This forces the two bundles to be of rank $n$,
and $\Cal L$ induces isomorphisms $R\cong L^*\otimes (TM/H)$ and $L\cong
R^*\otimes (TM/H)$.

In view of the description of $G_0$ above, the following result is a
special case of general prolongation procedures \cite{Tan,Mor,C-S},
see \cite{Takeuchi} and section 4.1 of \cite{C-tw} for more
information on this specific case.
\begin{thm*}
Let $H=L\oplus R$ be a Lagrangean contact structure on a manifold $M$
of dimension $2n+1$. Then there exists a principal $P$--bundle $p:\Cal
G\to M$ endowed with a Cartan connection $\om\in\Om^1(\Cal G,\g)$ such
that $L=Tp(\om^{-1}(\fg_{-1}^L\oplus\frak p))$ and
$R=Tp(\om^{-1}(\fg_{-1}^R\oplus\frak p))$. The pair $(\Cal G,\om)$ is uniquely
determined up to isomorphism provided that one in addition requires
the curvature of $\om$ to satisfy a normalization condition 
discussed in \ref{3.6}.  
\end{thm*}

Similarly, for any contact grading of a simple Lie algebra $\frak g$
and a choice of a Lie group $G$ with Lie algebra $\frak g$, one
defines a subgroup $P\subset G$ with Lie algebra
$\fg_0\oplus\fg_1\oplus\fg_2$. One then obtains an equivalence of
categories between regular normal parabolic geometries of type $(G,P)$
and underlying geometric structures, which in particular include a
contact structure. 

The second case of such structures we will be concerned with in this
paper, is partially integrable almost CR structures of hypersurface
type, see section \ref{5}. 

\subsection{Chains}\label{2.2}
Let $(p:\Cal G\to M,\om)$ be the canonical Cartan geometry determined
by a parabolic contact structure. Then one obtains an isomorphism
$TM\cong \Cal G\x_P(\frak g/\frak p)$ such that $H\subset TM$
corresponds to $(\frak g_{-1}\oplus\frak p)/\frak p\subset\frak g/\frak
p$. Of course, we may identify $\frak g_{-2}\oplus\frak g_{-1}$ as a
vector space with $\frak g/\frak p$ and use this to carry over the
natural $P$--action to $\g_{-2}\oplus\g_{-1}$. Let $Q\subset P$ be the
stabilizer of the line $\g_{-2}$ under this action. By definition,
this is a closed subgroup of $P$. Let us denote by $G_0\subset P$ the
closed subgroup consisting of all elements whose adjoint action
respects the grading of $\frak g$. Then $G_0$ has Lie algebra $\fg_0$
and by Proposition 2.10 of \cite{C-S}, any element $g\in P$ can be
uniquely written in the form $g_0\exp(Z_1)\exp(Z_2)$ for $g_0\in G_0$,
$Z_1\in\fg_1$, and $Z_2\in\fg_2$.

\begin{lemma*}
(1) An element $g=g_0\exp(Z_1)\exp(Z_2)\in P$ lies in the subgroup
  $Q\subset P$ if and only if $Z_1=0$. In particular, $\frak q=\frak
  g_0\oplus\fg_2$ and for $g\in Q$ we have
  $\Ad(g)(\fg_{-2})\subset\fg_{-2}\oplus\frak q$.

\noindent 
(2) Let $(p:\Cal G\to M,\om)$ be the canonical Cartan geometry
determined by a parabolic contact structure.  Let $x\in M$ be a point
and $\xi\in T_xM\setminus H_x$ a tangent vector transverse to the
contact subbundle.
  
Then there is a point $u\in p^{-1}(x)\subset\Cal G$ and a unique lift
$\tilde\xi\in T_u\Cal G$ of $\xi$ such that
$\om(u)(\tilde\xi)\in\g_{-2}$. The point $u$ is unique up to the
principal right action of an element $g\in Q\subset P$.
\end{lemma*}
\begin{proof}
  (1) We first observe that for a nonzero element $X\in\frak g_{-2}$,
  the map $Z\mapsto [Z,X]$ is a bijection $\fg_1\to\fg_{-1}$. This is
  easy to verify directly for the examples discussed in \ref{2.1} and
  \ref{5.1}. For general contact gradings it follows from the fact
  that $[\g_{-2},\g_2]$ consists of all multiples of the grading
  element, see section 4.2 of \cite{Yam}.
  
  By definition, $g\in Q$ if and only if
  $\Ad(g)(\fg_{-2})\subset\fg_{-2}\oplus\frak p$. Now from the
  expression $g^{-1}=\exp(-Z_2)\exp(-Z_1)g_0^{-1}$ one immediately
  concludes that $\Ad(g^{-1})(X)$ is congruent to
  $-[Z_1,X]\in\fg_{-1}$ modulo $\fg_{-2}\oplus\frak p$. Hence we see
  that $g\in Q$ if and only if $Z_1=0$, and the rest of (1) evidently
  follows.

\noindent
(2) Choose any point $v\in p^{-1}(x)$. Since the vertical bundle of
$\Cal G\to M$ equals $\om^{-1}(\frak p)$, there is a unique lift
$\eta\in T_v\Cal G$ of $\xi$ such that
$\om(v)(\eta)\in\g_{-2}\oplus\g_{-1}$. The assumption that $\xi$
is transverse to $H_x$ means that $\om(v)(\eta)\notin \g_{-1}$. For an
element $g\in P$ we can consider $v\cdot g$ and $T_vr^g\cdot\eta\in
T_{v\cdot g}\Cal G$, where $v\cdot g=r^g(v)$ denotes the principal
right action of $g$ on $v$.  Evidently, $T_vr^g\cdot\eta$ is again a
lift of $\xi$ and equivariancy of $\om$ implies that $\om(v\cdot
g)(T_vr^g\cdot\eta)=\Ad(g^{-1})(\om(v)(\eta))$.
  
Writing $\om(v)(\eta)=X_{-2}+X_{-1}$ we have $X_{-2}\neq 0$, so
from above we see that there is an element $Z\in\frak g_1$ such that
$[Z,X_{-2}]=X_{-1}$. Putting $g=\exp(Z)\in P$ we conclude that
$\om(v\cdot g)(T_vr^g\cdot\eta)\in\g_{-2}\oplus\frak p$. 
Hence putting $u=v\cdot g$ and subtracting
an appropriate vertical vector from $T_vr^g\cdot\eta$, we have found a
couple $(u,\tilde\xi)$ as required.

Any other choice of a preimage of $x$ has the form $u\cdot g$ for some
$g\in P$. Any lift of $\xi$ in $T_{u\cdot g}\Cal G$ is of the form
$T_ur^g\cdot\tilde\xi+\zeta$ for some vertical vector $\zeta$.
Clearly, there is a choice for $\zeta$ such that
$\om(T_ur^g\cdot\tilde\xi+\zeta)\in\g_{-2}$ if and only if
$\om(T_ur^g\cdot\tilde\xi)\in\g_{-2}\oplus\frak p$ and equivariancy of
$\om$ implies that this is equivalent to $g\in Q$.
\end{proof}

This lemma immediately leads us to chains: Fix a nonzero element
$X\in\fg_{-2}$. For a point $x\in M$ and a line $\ell$ in $T_xM$ which
is transverse to $H_x$, we can find a point $u\in\Cal G$ such that
$T_up\cdot\om_u^{-1}(X)\in\ell$. Denoting by $\tilde X$ the ``constant
vector field'' $\om^{-1}(X)$ we can consider the flow of $\tilde X$
through $u$ and project it onto $M$ to obtain a (locally defined)
smooth curve through $x$  whose tangent space at $x$ is $\ell$. In
section 4 of \cite{C-S-Z} it has been shown that, as an unparametrized
curve, this is uniquely determined by $x$ and $\ell$, and it comes
with a distinguished projective family of parametrizations. 

The lemma also leads us to a nice description of the space of all
transverse directions: For a point $u\in\Cal G$, we obtain a line in
$T_{p(u)}M$ which is transverse to $H_{p(u)}$, namely
$T_p(\om^{-1}_u(\fg_{-2}))$. This defines a smooth map $\Cal G\to\Cal
PTM$, where $\Cal PTM$ denotes the projectivized tangent bundle of
$M$. Since $P$ acts freely on $\Cal G$ so does $Q$ and hence $\Cal
G/Q$ is a smooth manifold. By the lemma, we obtain a
diffeomorphism from $\Cal G/Q$ to the open subset $\Cal
P_0TM\subset\Cal PTM$ formed by all lines which are transverse to the
contact distribution $H$.

\subsection{Path geometries}\label{2.3}
Classically, path geometries are associated to certain families of
unparametrized curves in a smooth manifold. Suppose that in a manifold
$Z$ we have a smooth family of curves such that through each point of
$Z$ there is exactly one curve in each direction. Let $\Cal PTZ$ be
the projectivized tangent bundle of $Z$, i.e.~the space of all lines
through the origin in tangent spaces of $Z$. Given a line $\ell$ in
$T_xZ$, we can choose the unique curve in the family which goes
through $x$ in direction $\ell$. Choosing a local regular
parametrization $c:I\to Z$ of this curve we obtain a lift $\tilde
c:I\to\Cal PTZ$ by defining $\tilde c(t)$ to be the line in $T_{c(t)}Z$
generated by $c'(t)$. Choosing a different regular parametrization, we
just obtain a reparametrization of $\tilde c$, so the submanifold
$\tilde c(I)\subset\Cal PTZ$ is independent of all choices. These
curves foliate $\Cal PTZ$, and their tangent spaces give rise to a
line subbundle $E\subset T\Cal PTZ$. 

This subbundle has a special property: Similarly to the tautological
line bundle on a projective space, a projectivized tangent bundle
carries a tautological subbundle $\Xi\subset T\Cal PTZ$ of rank
$\dim(Z)$. By definition, given a line $\ell\subset T_zZ$, a tangent
vector $\xi\in T_\ell\Cal PTZ$ lies in $\Xi_\ell$ if and only if its
image under the tangent map of the projection $\Cal PTZ\to Z$ lies in
the line $\ell$. By construction, the line subbundle $E$ associated to
a family of curves as above always is contained in $\Xi$ and is
transverse to the vertical subbundle $V$ of $\Cal PTZ\to Z$. Hence we
see that $\Xi=E\oplus V$.

Conversely, having given a decomposition $\Xi=E\oplus V$ of the
tautological bundle, we can project the leaves of the foliation of
$\Cal PTZ$ defined by $E$ to the manifold $Z$ to obtain a smooth
family of curves in $Z$ with exactly one curve through each point in
each direction. Hence one may use the decomposition $\Xi=E\oplus V$ as
an alternative definition of such a family of curves, and this
decomposition is usually referred to as a \textit{path geometry} on
$Z$. It is easy to verify that the Lie bracket of vector fields
induces an isomorphism $E\otimes V\to T\Cal PTZ/\Xi$. 

It turns out that path geometries also admit an equivalent description
as regular normal parabolic geometries. Putting $m:=\dim(Z)-1$ we
consider the Lie algebra $\tg:=\frak {sl}(m+2,\Bbb R)$ with the
$|2|$--grading obtained by a block decomposition 
$$
\begin{pmatrix}
\tg_0&\tg^E_1&\tg_2\\\tg^E_{-1}&\tg_0&\tg^V_1\\\tg_{-2}&\tg^V_{-1}&\tg_0
\end{pmatrix}.
$$
as in \ref{2.1}, but this time with blocks of size $1$, $1$, and
$m$. Hence $\tg_{\pm 1}^E$ has dimension $1$ while $\tg_{\pm 1}^V$ and
$\tg_{\pm 2}$ are all $m$--dimensional. Put $\tilde G:=PGL(m+2,\Bbb
R)$ and let $\tilde G_0\subset\tilde P\subset\tilde G$ be the
subgroups formed by matrices which are block diagonal respectively
block upper triangular with block sizes $1$, $1$, and $m$. Then
$\tilde G_0$ and $\tilde P$ have Lie algebras $\tg_0$ respectively
$\tp:=\tg_0\oplus\tg_1\oplus\tg_2$, where $\tg_1=\tg_1^E\oplus\tg_1^V$. 

The adjoint action identifies $\tilde G_0$ with the group of
automorphisms of the graded Lie algebra $\tg_{-2}\oplus\tg_{-1}$ which
in addition preserve the decomposition $\tg_{-1}=\tg^E_{-1}\oplus
\tg^V_{-1}$.  Hence the following result is a special case of the
general prolongation procedures \cite{Tan,Mor,C-S}, see section 4.7 of
\cite{C-tw} for this specific case.
\begin{thm*}
  Let $\tilde Z$ be a smooth manifold of dimension $2m+1$ endowed with
  transversal subbundles $E$ and $V$ in $T\tilde Z$ of rank $1$ and
  $m$, respectively, and put $\Xi:=E\oplus V\subset T\tilde Z$. Suppose
  that the Lie bracket of two sections of $V$ is a section of $\Xi$
  and that the tensorial map $E\otimes V\to T\tilde Z/\Xi$ induced by
  the Lie bracket of vector fields is an isomorphism.
  
  Then there exists a principal bundle $\tilde p:\tcg\to\tilde Z$ with
  structure group $\tilde P$ endowed with a Cartan connection
  $\tilde\om\in\Om^1(\tcg,\tg)$ such that $E=T\tilde
  p(\tilde\om^{-1}(\tg_{-1}^E\oplus\tp))$ and $V=T\tilde
  p(\tilde\om^{-1}(\tg_{-1}^V\oplus\tp))$. The pair $(\tcg,\tilde\om)$
  is uniquely determined up to isomorphism provided that $\tilde\om$
  is required to satisfy a normalization condition discussed in \ref{3.6}. 
\end{thm*}
In particular, a family of paths on $Z$ as before gives rise to a
Cartan geometry on $\Cal PTZ$. This immediately generalizes to the
case of an open subset of $\Cal PTZ$, i.e.~the case where paths are
only given through each point in an open set of directions. 

It turns out that for $m\neq 2$, the assumptions of the theorem
already imply that the subbundle $V\subset T\tilde Z$ is involutive.
Then $\tilde Z$ is automatically locally diffeomorphic to a
projectivized tangent bundle in such a way that $V$ is mapped to the
vertical subbundle and $\Xi$ to the tautological subbundle. Hence for
$m\neq 2$, the geometries discussed in the theorem are locally
isomorphic to path geometries. 

\subsection{The path geometry of chains}\label{2.4}
From \ref{2.2} we see that for a manifold $M$ endowed with a parabolic
contact structure the chains give rise to a path geometry on the open
subset $\tilde M:=\Cal P_0TM$ of the projectivized tangent bundle of
$M$. We can easily describe the corresponding configuration of bundles
explicitly: Denoting by $(p:\Cal G\to M,\om)$ the Cartan geometry
induced by the parabolic contact structure, we know from \ref{2.2}
that $\tilde M=\Cal G/Q$, where $Q\subset P$ denotes the stabilizer of
the line in $\fg/\frak p$ corresponding to
$\fg_{-2}\subset\fg_{-2}\oplus\fg_{-1}$. In particular, $\Cal G$ is a
$Q$--principal bundle over $\tilde M$ and $\om$ is a Cartan connection
on $\Cal G\to\tilde M$. This implies that $T\tilde M=\Cal G\x_Q\frak
g/\frak q$, and the tangent map to the projection $\pi:\tilde M\to M$
corresponds to the obvious projection $\frak g/\frak q\to\frak g/\frak
p$. In particular, the vertical bundle $V=\ker(T\pi)$ corresponds to
$\frak p/\frak q\subset\frak g/\frak q$. From the construction of the
isomorphism $\Cal G/Q\to\tilde M$ in \ref{2.2}, it is evident that the
tautological bundle $\Xi$ corresponds to $(\frak g_{-2}\oplus\frak
p)/\frak q$.  By part (1) of Lemma \ref{2.2}, the subspace $(\frak
g_{-2}\oplus\frak q)/\frak q\subset\frak g/\frak q$ is $Q$--invariant,
thus it gives rise to a line subbundle $E$ in $\Xi$, which is
complementary to $V$. By construction, this exactly describes the path
geometry determined by the chains.

If $\dim(M)=2n+1$, then the dimension of $\tilde M$ is $4n+1$. Put
$\tilde G:=PGL(2n+2,\Bbb R)$ and let $\tilde P\subset\tilde G$ be the
subgroup described in \ref{2.3}. Then by Theorem \ref{2.3} the path
geometry on $\tilde M$ gives rise to a canonical principal bundle
$\tcg\to\tilde M$ with structure group $\tilde P$ endowed with a
canonical normal Cartan connection $\tilde\om\in\Om^1(\tcg,\tg)$. The
main question now is whether there is a direct relation between the
Cartan geometries $(\Cal G\to\tilde M,\om)$ and $(\tcg\to
\tilde M,\tilde\om)$. 

The only reasonable way to relate these two Cartan geometries is to
consider a morphism $j:\Cal G\to\tcg$ of principal bundles and compare
the pull--back $j^*\tilde\om$ to $\om$. This means that $j$ is
equivariant, so we first have to choose a group homomorphism
$i:Q\to\tilde P$ and require that $j(u\cdot g)=j(u)\cdot i(g)$ for all
$g\in Q$. Having chosen $i$ and $j$, we have
$j^*\tilde\om\in\Om^1(\Cal G,\tg)$ and the only way to directly relate
this to $\om\in\Om^1(\Cal G,\frak g)$ is to have
$j^*\tilde\om=\al\o\om$ for some linear map $\al:\frak g\to\tg$. If we
have such a relation, then we can immediately recover $\tcg$ from
$\Cal G$: Consider the map $\Ph:\Cal G\x\tilde P\to\tcg$ defined by
$\Ph(u,\tilde g):=j(u)\cdot\tilde g$. Equivariancy of $j$ immediately
implies that $\Ph(u\cdot g,\tilde g)=\Ph(u,i(g)\tilde g)$, so $\Ph$
descends to a bundle map $\Cal G\x_Q\tilde P\to\tcg$, where the left
action of $Q$ on $\tilde P$ is defined via $i$. This is immediately
seen to be an isomorphism of principal bundles, so $\tcg$ is obtained
from $\Cal G$ by an extension of structure group. Under this
isomorphism, the given morphism $j:\Cal G\to\tcg$ corresponds to the
natural inclusion $\Cal G\to\Cal G\x_Q\tilde P$ induced by $u\mapsto
(u,e)$.

\section{Induced Cartan connections}\label{3}
In this section, we study the problem of extending Cartan connections.
We derive the basic results in the setting of general Cartan
geometries, and then specialize to the case of parabolic contact
structures and, in particular, Lagrangean contact structures. Some of
the developments in \ref{3.1} and \ref{3.3} below are closely related
to \cite{Kobayashi, Wang}.

\subsection{Extension functors for Cartan geometries}\label{3.1}
Motivated by the last observations in \ref{2.4}, let us consider the
following problem: Suppose we have given Lie groups $G$ and $\tilde G$
with Lie algebras $\frak g$ and $\tg$, closed subgroups $Q\subset G$
and $\tilde P\subset\tilde G$, a homomorphism $i:Q\to\tilde P$ and a
linear map $\al:\frak g\to\tg$. We will assume throughout that $i$ is
infinitesimally injective, i.e.~$i':\frak q\to\tp$ is injective. 

Given a Cartan geometry $(p:\Cal G\to N,\om)$ of type $(G,Q)$, we put
$\tcg:=\Cal G\x_Q\tilde P$ and denote by $j:\Cal G\to\tcg$ the
canonical map. Since $i$ is infinitesimally injective, this is an
immersion, i.e.~$T_uj$ is injective for all $u\in\Cal G$. We want to
understand whether there is a Cartan connection
$\tilde\om\in\Om^1(\tcg,\tg)$ such that $j^*\tilde\om=\al\o\om$, and
if so, whether $\tilde\om$ is uniquely determined.

\begin{prop*}
  There is a Cartan connection $\tilde\om$ on $\tcg$ such that
  $j^*\tilde\om=\al\o\om$ if and only if the pair $(i,\al)$ satisfies
  the following conditions:
\begin{enumerate}
\item $\al\o\Ad(g)=\Ad(i(g))\o\al$ for all $g\in Q$.
\item On the subspace $\frak q\subset\frak g$, the map $\al$ restricts
  to the derivative $i'$ of $i:Q\to\tilde P$.
\item The map $\underline{\al}:\frak g/\frak q\to\tg/\tp$ induced by
  $\al$ is a linear isomorphism.
\end{enumerate}
If these conditions are satisfied, then $\tilde\om$ is uniquely
determined. 
\end{prop*}
\begin{proof}
  Let us first assume that there is a Cartan connection $\tilde\om$ on
  $\tcg$ such that $j^*\tilde\om=\al\o\om$. For $u\in\Cal G$, the
  tangent space $T_{j(u)}\tcg$ is spanned by $T_uj(T_u\Cal G)$ and the
  vertical subspace $V_{j(u)}\tcg$. The behavior of $\tilde\om$ on the
  first subspace is determined by the fact that
  $j^*\tilde\om=\al\o\om$, while on the second subspace $\tilde\om$
  has to reproduce the generators of fundamental vector fields. Hence
  the restriction of $\tilde\om$ to $j(\Cal G)$ is determined by the
  fact that $j^*\tilde\om=\al\o\om$. By definition of $\tcg$, any
  point $\tilde u\in\tcg$ can be written as $j(u)\cdot\tilde g$ for
  some $u\in\Cal G$ and some $\tilde g\in\tilde P$, so uniqueness of
  $\tilde\om$ follows from equivariancy. 
  
  Still assuming that $\tilde\om$ exists, condition (1) follows from
  equivariancy of $j$, $\om$, and $\tilde\om$. Equivariancy of $j$
  also implies that for $A\in\frak q$ and the corresponding
  fundamental vector field $\ze_A$ we get $Tj\o\ze_A=\ze_{i'(A)}$.
  Thus condition (2) follows from the fact that both $\om$ and
  $\tilde\om$ reproduce the generators of fundamental vector fields.
  Let $p:\Cal G\to N$ and $\tilde p:\tcg\to N$ be the bundle
  projections, so $\tilde p\o j=p$. For $\xi\in T_u\Cal G$ we have
  $\al(\om(\xi))=\tilde\om(T_uj\cdot\xi)$, so if this lies in $\tp$
  then $T_uj\cdot\xi$ is vertical. But then $\xi$ is vertical and
  hence $\om(\xi)\in\frak q$. Therefore, the map $\underline{\al}$ is
  injective, and since both $\Cal G$ and $\tcg$ admit a Cartan
  connection, we must have $\dim(\frak g/\frak
  q)=\dim(N)=\dim(\tg/\tp)$, so (3) follows.
  
  Conversely, suppose that (1)--(3) are satisfied for $(i,\al)$ and
  $\om$ is given. For $\tilde u\in\tcg$ and $\tilde\xi\in T_{\tilde
    u}\tcg$ we can find elements $u\in\Cal G$, $\xi\in T_u\Cal G$,
  $A\in\tp$, and $\tilde g\in\tilde P$ such that $\tilde
  u=j(u)\cdot\tilde g$ and $\tilde\xi=Tr^{\tilde g}\cdot
  (Tj\cdot\xi+\ze_A)$.  Then we define
  $\tilde\om(\tilde\xi):=\Ad(\tilde g)^{-1}(\al(\om(\xi))+A)$. Using
  properties (1) and (2) one verifies that this is independent of all
  choices. By (3), it defines a linear isomorphism $T_{\tilde
    u}\tcg\to\tg$, and the remaining properties of a Cartan connection
  are easily verified directly.
\end{proof}

Any pair $(i,\al)$ which satisfies the properties (1)--(3) of the
proposition gives rise to an extension functor from Cartan geometries
of type $(G,Q)$ to Cartan geometries of type $(\tilde G,\tilde P)$:
Starting from a geometry $(p:\Cal G\to N,\om)$ of type $(G,Q)$, one
puts $\tcg:=\Cal G\x_Q\tilde P$ (with $Q$ acting on $\tilde P$ via
$i$) and defines $\tilde\om\in\Om^1(\tcg,\tg)$ to be the unique Cartan
connection on $\tcg$ such that $j^*\tilde\om=\al\o\om$, where $j:\Cal
G\to\tcg$ is the canonical map. For a morphism $\ph:\Cal
G_1\to\Cal G_2$ between geometries of type $(G,Q)$, we can consider the
principal bundle map $\Ph:\tcg_1\to\tcg_2$ induced by $\ph\x\id_{\tilde
  P}$. By construction, this satisfies $\Ph\o j_1=j_2\o\ph$ and we
obtain 
$$
j_1^*\Ph^*\tilde\om_2=\ph^*j_2^*\tilde\om_2=\ph^*(\al\o\om_2)=
\al\o\ph^*\om_2=\al\o\om_1. 
$$
But $\tilde\om_1$ is the unique Cartan connection whose pull--back
along $j_1$ coincides with $\al\o\om_1$, which implies that
$\Ph^*\tilde\om_2=\tilde\om_1$, and hence $\Ph$ is a morphism of
Cartan geometries of type $(\tilde G,\tilde P)$.

There is a simple notion of equivalence for pairs $(i,\al)$: We call
$(i,\al)$ and $(\hat i,\hat\al)$ equivalent and write
$(i,\al)\sim(\hat i,\hat\al)$ if and only if there is an element
$\tilde g\in\tilde P$ such that $\hat i(g)=\tilde g^{-1}i(g)\tilde g$
and $\hat\al=\Ad(\tilde g^{-1})\o\al$. Notice that if $(i,\al)$
satisfies conditions (1)--(3) of the proposition, then so does any
equivalent pair. In order to distinguish between different extension
functors, for a geometry $(p:\Cal G\to M,\om)$ of type $(G,Q)$ we will
often denote the geometry of type $(\tilde G,\tilde P)$ obtained using
$(i,\al)$ by $(\Cal G\x_i\tilde P,\tilde\om_\al)$.

\begin{lemma*}
  Let $(i,\al)$ and $(\hat i,\hat\al)$ be equivalent pairs satisfying
  conditions (1)--(3) of the proposition. Then the resulting extension
  functors for Cartan geometries are naturally isomorphic.
\end{lemma*}
\begin{proof}
  By assumption, there is an element $\tilde g\in\tilde P$ such that
  $\hat i(g)=\tilde g^{-1}i(g)\tilde g$ and $\hat\al=\Ad(\tilde
  g^{-1})\o\al$. Let $j:\Cal G\to\Cal G\x_i\tilde P$ and $\hat j:\Cal
  G\to\Cal G\x_{\hat i}\tilde P$ be the natural inclusions, and
  consider the map $r^{\tilde g}\o j:\Cal G\to\Cal G\x_i\tilde
  P$. Evidently, we have $j(u\cdot g)\cdot\tilde g=j(u)\cdot\tilde
  g\cdot\hat i(g)$.  Hence, by the last observation in \ref{2.4}, we
  obtain an isomorphism $\Ps:\Cal G\x_{\hat i}\tilde P\to\Cal
  G\x_i\tilde P$ such that $\Ps\o\hat j=r^{\tilde g}\o j$. Now we
  compute
$$
  \hat j^*\Ps^*\tilde\om_\al=j^*(r^{\tilde g})^*\tilde\om_\al=
  \Ad(\tilde g^{-1})\o j^*\tilde\om_\al=\hat\al\o\om.
$$
  By uniqueness, $\Ps^*\tilde\om_\al=\tilde\om_{\hat\al}$, so $\Ps$
  is a morphism of Cartan geometries. It is clear from the
  construction that this defines a natural transformation between the
  two extension functors and an inverse can be constructed in the same
  way using $\tilde g^{-1}$ rather than $\tilde g$.
\end{proof}

\subsection{The relation to the Fefferman construction}\label{3.2}
There is a simple source of pairs $(i,\al)$ which satisfy conditions
(1)--(3) of Proposition \ref{3.1}: Suppose that $\ph:G\to\tilde G$ is
an infinitesimally injective homomorphism of Lie groups such that
$\ph(Q)\subset\tilde P$. Then $i:=\ph|_Q:Q\to\tilde P$ is an
infinitesimally injective homomorphism and $\al:=\ph':\frak g\to\tg$
is a Lie algebra homomorphism. Then condition (2) of Proposition
\ref{3.1} is satisfied by construction, while condition (1) easily
follows from differentiating the equation
$\ph(ghg^{-1})=\ph(g)\ph(h)\ph(g)^{-1}$. Hence the only nontrivial
condition is (3). Note that if $(i,\al)$ is obtained from $\ph$ in
this way, than any pair equivalent to $(i,\al)$ is obtained in the
same way from the map
$g\mapsto \tilde g\ph(g)\tilde g^{-1}$ for some $\tilde
g\in\tilde G$. The main feature of such pairs is that $\al$ is a
homomorphism of Lie algebras.

In this setting, one may actually go one step further: Suppose we have
fixed an infinitesimally injective $\ph:G\to\tilde G$ and a closed
subgroup $\tilde P\subset\tilde G$. Then we put $Q:=\ph^{-1}(\tilde
P)\subset G$ to obtain a pair $(i:=\ph|_Q,\al:=\ph')$ and hence an
extension functor from Cartan geometries of type $(G,Q)$ to geometries
of type $(\tilde G,\tilde P)$. For a closed subgroup $P\subset G$ with
$Q\subset P$, one gets a functor from geometries of type $(G,P)$ to
geometries of type $(G,Q)$ as described in \ref{2.2}: Given a geometry
$(p:\Cal G\to M,\om)$ of type $(G,P)$, one defines $\tilde M:=\Cal
G/Q=\Cal G\x_P(P/Q)$ and $(\Cal G\to\tilde M,\om)$ is a geometry of
type $(G,Q)$. Combining with the above, one gets a functor from
geometries of type $(G,P)$ to geometries of type $(\tilde G,\tilde
P)$. 

The most important example of this is the Cartan geometry
interpretation of Fefferman's construction of a Lorentzian conformal
structure on the total space of a certain circle bundle over a CR
manifold, see \cite{Fefferman}. In this case $G=SU(n+1,1)$, $\tilde
G=SO(2n+2,2)$, and $\ph$ is the evident inclusion. Putting $\tilde P$
the stabilizer of a real null line $\ell\subset\Bbb R^{2n+4}$ in
$\tilde G$, the group $Q=G\cap\tilde P$ is the stabilizer of $\ell$ in
$G$. Evidently, this is contained in the stabilizer $P\subset G$ of
the complex null line spanned by $\ell$, and $P/Q\cong\Bbb RP^1\cong
S^1$. Hence the above procedure defines a functor, which to a
parabolic geometry of type $(G,P)$ on $M$ associates a parabolic
geometry of type $(\tilde G,\tilde P)$ on the total space $\tilde M$
of a circle bundle over $M$. More details about this can be found in
\cite{Cap:Fefferman}. 

\subsection{The effect on curvature}\label{3.3}
We next discuss the effect of extension functors of the type discussed
in \ref{3.1} on the curvature of Cartan geometries. This will show
specific features of the special case discussed in \ref{3.2}. 

For a Cartan connection $\om$ on a principal $P$--bundle $\Cal G\to M$
with values in $\frak g$, one initially defines the curvature
$K\in\Om^2(\Cal G,\frak g)$ by
$K(\xi,\eta):=d\om(\xi,\eta)+[\om(\xi),\om(\eta)]$. This measures the
amount to which the Maurer--Cartan equation fails to hold.  The
defining properties of a Cartan connection immediately imply that $K$
is horizontal and $P$--equivariant. In particular, $K(\xi,\eta)=0$ for
all $\eta$ provided that $\xi$ is vertical or, equivalently, that
$\om(\xi)\in\frak p$.

Using the trivialization of $T\Cal G$ provided by $\om$, one can pass
to the curvature function $\ka:\Cal G\to L(\La^2(\frak g/\frak
p),\frak g)$, which is characterized by
$$
\ka(u)(X+\frak p,Y+\frak p):=K(u)(\om^{-1}(X),\om^{-1}(Y)).
$$
This is well defined by horizontality of $K$, and equivariancy of $K$
easily implies that $\ka$ is equivariant for the natural $P$--action on the
space $L(\La^2(\frak g/\frak p),\frak g)$, which is induced from the
adjoint action on all copies of $\frak g$. 

Using the setting of \ref{3.1}, suppose that $(i:Q\to\tilde
P,\al:\frak g\to\tg)$ is a pair satisfying the conditions (1)--(3) of
Proposition \ref{3.1}. Consider the map $\frak g\x\frak g\to \tg$
defined by $(X,Y)\mapsto [\al(X),\al(Y)]_{\tg}-\al([X,Y]_{\frak g})$,
which measures the deviation from $\al$ being a homomorphism of Lie
algebras. This map is evidently skew symmetric. By condition (1),
$\al\o\Ad(g)=\Ad(i(g))\o\al$ for all $g\in Q$, which infinitesimally
implies that $\al\o\ad(X)=\ad(i'(X))\o\al$ for all $X\in\frak q$, and
by condition (2) we have $i'(X)=\al(X)$ in this case. Hence this map
vanishes if one of the entries is from $\frak q\subset\frak g$, and we
obtain a well defined linear map $\La^2(\frak g/\frak q)\to \tg$. By
condition (3), $\al$ induces a linear isomorphism
$\underline{\al}:\frak g/\frak q\to\tg/\tp$, and we conclude that we
obtain a well defined map $\Ps_{\al}:\La^2(\tg/\tp)\to\tg$
by putting 
\begin{equation*}
\Ps_{\al}(\tilde X+\tp,\tilde Y+\tp)=[\al(X),\al(Y)]-\al([X,Y]), 
\end{equation*}
where $\al(X)+\tp=\tilde X+\tp$ and $\al(Y)+\tp=\tilde Y+\tp$. 

\begin{prop*}
  Let $(i,\al)$ be a pair satisfying conditions (1)--(3) of
  Proposition \ref{3.1}. Let $(p:\Cal G\to N,\om)$ be a Cartan
  geometry of type $(G,Q)$, let $(\Cal G\x_i\tilde P,\tilde\om_{\al})$ be
  the geometry of type $(\tilde G,\tilde P)$ obtained using the
  extension functor associated to $(i,\al)$, and let $j:\Cal G\to\Cal
  G\x_i\tilde P$ be the natural map.

Then the curvature functions $\ka$ and $\tilde\ka$ of the two
geometries satisfy 
$$
\tilde\ka(j(u))(\tilde X,\tilde
Y)=\al(\ka(u)(\underline{\al}^{-1}(\tilde
X),\underline{\al}^{-1}(\tilde Y)))+\Ps_{\al}(\tilde X,\tilde Y), 
$$ 
for any $\tilde X,\tilde Y\in\tg/\tp$,
and this completely determines $\tilde\ka$. 

In particular, if $\om$ is flat, then $\tilde\om$ is flat if and only
if $\al$ is a homomorphism of Lie algebras. 
\end{prop*}
\begin{proof}
By definition, $j^*\tilde\om_\al=\al\o\om$, and hence
$j^*d\tilde\om_{\al}=\al\o d\om$. This immediately implies that for the
curvatures $K$ and $\tilde K$ and $\xi,\eta\in\frak X(\Cal G)$ we get
$$
  \tilde K(j(u))(Tj\cdot\xi,Tj\cdot\eta)=\al(d\om(u)(\xi,\eta))+
[\al(\om(u)(\xi)),\al(\om(u)(\eta))].
$$
On the other hand, we get
$$
\al(K(u)(\xi,\eta))=\al(d\om(u)(\xi,\eta))+\al([\om(u)(\xi),\om(u)(\eta)]).
$$
Now the formula for $\tilde\ka(j(u))$ follows immediately from the
definition of the curvature functions. Since $\tilde\ka$ is $\tilde
P$--equivariant, it is completely determined by its restriction to
$j(\Cal G)$. The final claim follows directly, since $\Ps_\al$ vanishes
if and only if $\al$ is a homomorphism of Lie algebras. 
\end{proof}

\subsection{Uniqueness}\label{3.4}
A crucial fact for the further development is that, passing from
parabolic contact structures to the associated path geometries of
chains, there is actually no freedom in the choice of the pair
$(i,\al)$ up to equivalence as introduced in \ref{3.1} above. This
result certainly is valid in a more general setting but it seems to be
difficult to give a nice formulation for conditions one has to assume. 

Therefore we return to the setting of section \ref{2}, i.e.~$G$ is
semisimple, $P\subset G$ is obtained from a contact grading, $Q$ is
the subgroup described in \ref{2.2}, and $\tilde G$ and $\tilde P$
correspond to path geometries in the appropriate dimension as in
\ref{2.3}. In this setting we can now prove:

\begin{thm*}
  Let $(i,\al)$ and $(\hat i,\hat\al)$ be pairs satisfying conditions
  (1)--(3) of Proposition \ref{3.1}. Suppose that there is a Cartan
  geometry $(p:\Cal G\to M,\om)$ of type $(G,Q)$ such that there is an
  isomorphism between the geometries of type $(\tilde G,\tilde P)$
  obtained using $(i,\al)$ and $(\hat i,\hat\al)$, which covers the
  identity on $M$. Then $(i,\al)$ and $(\hat i,\hat\al)$ are
  equivalent. 
\end{thm*}
\begin{proof}
  Using the notation of the proof of Lemma \ref{3.1}, suppose that we
  have an isomorphism $\Ps:\Cal G\x_i\tilde P\to\Cal G\x_{\hat
    i}\tilde P$ of principal bundles which covers the identity on $M$
  and has the property that
  $\Ps^*\tilde\om_{\hat\al}=\tilde\om_{\al}$. Let us denote by $j$ and
  $\hat j$ the natural inclusions of $\Cal G$ into the two extended
  bundles. Since $\Ps$ covers the identity on $M$, there must be a
  smooth function $\ph:\Cal G\to \tilde P$ such that $\Ps(j(u))=\hat
  j(u)\cdot\ph(u)$.
    
  By construction we have $j(u\cdot g)=j(u)\cdot i(g)$ and $\hat
  j(u\cdot g)=\hat j(u)\cdot\hat i(g)$, and using the fact that $\Ps$
  is $\tilde P$--equivariant we obtain $\hat i(g)=\ph(u)i(g)\ph(u\cdot
  g)^{-1}$. On the other hand, differentiating the equation
  $\Ps(j(u))=\hat j(u)\cdot\ph(u)$, we obtain
$$
(T\Ps\o Tj)\cdot\xi=(Tr^{\ph(u)}\o T\hat
j)\cdot\xi+\ze_{\delta\ph(u)(\xi)}(\Ps(j(u)))
$$
where $\delta\ph\in\Om^1(\Cal G,\tp)$ denotes the left logarithmic
derivative of $\ph:\Cal G\to \tilde P$. Applying $\tilde\om_{\hat\al}$
to the left hand side of this equation, we simply get
$$
(j^*\Ps^*\tilde\om_{\hat\al})(\xi)=(j^*\tilde\om_\al)(\xi)=\al(\om(\xi)).
$$
Applying $\tilde\om_{\hat\al}$ to the right hand side, we obtain 
\begin{multline*}
  (\hat j^*(r^{\ph(u)})^*\tilde\om_{\hat\al})(\xi)+\delta\ph(u)(\xi)=\\
\Ad(\ph(u)^{-1})((\hat j^*\tilde\om_{\hat\al})(\xi))+\delta\ph(u)(\xi)=\\
\Ad(\ph(u)^{-1})(\hat\al(\om(\xi)))+\delta\ph(u)(\xi), 
\end{multline*}
and we end up with the equation
\begin{equation*}
\al(\om(\xi))=\Ad(\ph(u)^{-1})(\hat\al(\om(\xi)))+\delta\ph(u)(\xi)
\tag{$*$}  
\end{equation*}
for all $\xi\in T\Cal G$. Together with the relation between $i$ and
$\hat i$ derived above, this shows that it suffices to show that
$\ph(u)$ is constant to prove that $(i,\al)\sim(\hat i,\hat\al)$. 

By construction, $\delta\ph(u)$ has values in $\tp$, so projecting
equation ($*$) to $\tg/\tp$ implies that
$$
\al(\om(\xi))+\tp=\Adb(\ph(u)^{-1})(\hat\al(\om(\xi))+\tp),
$$ 
for all $\xi\in T_u\Cal G$, where $\Adb$ is the action of $\tilde P$
on $\tg/\tp$ induced by the adjoint action. By property (3) from
Proposition \ref{3.1} this implies that
$\underline{\al}=\Adb(\ph(u)^{-1})\o\underline{\hat\al}$, so we see
that $\Adb(\ph(u)^{-1})$ must be independent of $u$. Hence we must
have $\ph(u)=\tilde g_1\ph_1(u)$ for some element $\tilde g_1\in\tilde
P$ and a smooth function $\ph_1:\Cal G\to \tilde P$ which has values
in the kernel of $\Adb$. As in \ref{2.2}, any element of $\tilde P$
can be uniquely written in the form $\tilde g_0\exp(\tilde
Z_1)\exp(\tilde Z_2)$ with $\tilde g_0\in\tilde G_0$ and $\tilde
Z_i\in\tg_i$, and such an element lies in the kernel of $\Adb$ if and
only if $\Ad(\tilde g_0)$ restricts to the identity on $\tg_-$ and
$\tilde Z_1=0$. Since $\tp_+$ is dual to $\tg_-$ and $\tg_0$ injects
into $L(\tg_-,\tg_-)$ the first condition implies that $\Ad(\tilde
g_0)=\id_{\tilde\g}$. Since $\tilde G=PGL(k,\Bbb R)$ for some $k$,
this implies that $\tilde g_0$ is the identity.

Hence $\ph_1$ has values in $\exp(\tg_2)$ and therefore $\delta\ph(u)$
has values in $\tg_2$. Projecting equation ($*$) to $\tg/\tg_2$, we
obtain
$$
\al(\om(\xi))+\tg_2=\Adb(\ph(u)^{-1})(\hat\al(\om(\xi))+\tg_2),
$$
where this time $\Adb$ denotes the natural action on $\tg/\tg_2$. But
by \cite[Lemma 3.2]{Yam} an element of $\tg_2$ vanishes provided that
all brackets with elements of $\tg_{-1}$ vanish, and this easily
implies that $\ph_1(u)$ is the identity and so $\ph$ is constant. 
\end{proof}

This result has immediate consequences on the problem of describing
the path geometry of chains associated to a parabolic contact
structure: If we start with the homogeneous model $G/P$ for a
parabolic contact geometry, the induced path geometry of chains is
defined on the homogeneous space $G/Q$. To obtain this by an extension
functor as described in \ref{3.1}, we need a homomorphism
$i:Q\to\tilde P$ and a linear map $\al:\frak g\to\tg$, where $(\tilde
G,\tilde P)$ gives rise to path geometries in the appropriate
dimension. The pair $(i,\al)$ has to satisfy conditions (1)--(3) of
Proposition \ref{3.1} in order to give rise to an extension functor.
The only additional condition is that the extended geometry
$(G\x_i\tilde P,\tilde\om_{\al})$ obtained from $(G\to G/Q,\om^{MC})$
is regular and normal. By Theorem \ref{2.3}, a regular normal
parabolic geometry of type $(\tilde G,\tilde P)$ is uniquely
determined by the underlying path geometry, which is encoded into
$(G\to G/Q,\om^{MC})$, see \ref{2.4}.

The theorem above then implies that $(i,\al)$ is uniquely determined
up to equivalence. In view of Lemma \ref{3.1}, the extension functor
obtained from $(i,\al)$ is (up to natural isomorphism) the only
extension functor of the type discussed in \ref{3.1} which produces
the right result for the homogeneous model (and hence for locally flat
geometries).

The final step is then to study under which conditions on a geometry
of type $(G,P)$, the extension functor associated to $(i,\al)$
produces a regular normal geometry of type $(\tilde G,\tilde P)$.

\subsection{}\label{3.5}
Let us return to the case of Lagrangean contact structures as
discussed in \ref{2.1}. By definition, we have $G=PGL(n+2,\Bbb R)$ and
$P\subset G$ is the subgroup of all matrices which are block upper
triangular with blocks of sizes $1$, $n$, and $1$. From part (1) of
Lemma \ref{2.2} one immediately concludes that $Q\subset P$ is the
subgroup formed by all matrices of the block form
%%% I have changed the letters, following 5.2:
$$
\begin{pmatrix}  p & 0 & s\\ 0 & R & 0\\ 0 & 0 & q\end{pmatrix},
$$
such that $|pq\det(R)|=1$. Since the corresponding manifolds have
dimension $2n+1$, the right group for the path geometry defined by the
chains is $\tilde G=PGL(2n+2,\Bbb R)$. The subgroup $\tilde
P\subset\tilde G$ is given by the classes of those matrices which are
block upper triangular with blocks of sizes $1$, $1$, $2n$. In the
sequel, we will always further split the last block into two blocks of
size $n$.

Consider the (well defined) smooth map $i:Q\to\tilde P$ and the linear map 
$\al:\fg\to\tg$ defined by
\begin{equation*}
  i\begin{pmatrix}p&0&s\\0&R&0\\0&0&q\end{pmatrix}:=
  \begin{pmatrix}\sgn(\tfrac{q}{p})\sqrt{|\tfrac{p}{q}|} &
  \sgn(\tfrac{q}{p})\tfrac{s}{p}\sqrt{|\tfrac{p}{q}|} & 0 & 0\\
  0 & \sqrt{|\tfrac{q}{p}|} & 0 & 0 \\
  0&0&q^{-1}\sqrt{|\tfrac{q}{p}|}R&0\\
  0&0&0&p\sqrt{|\tfrac{q}{p}|}(R^{-1})^t\end{pmatrix},
\end{equation*}
\begin{equation*}
 \al\begin{pmatrix}a&u&d\\x&B&v\\z&y&c\end{pmatrix}:
=\begin{pmatrix}\tfrac{a-c}{2}&d&\tfrac12 u&\tfrac12 v^t\\
z&\frac{c-a}{2}&\tfrac12 y&-\tfrac12 x^t\\
x&v&B-\frac{a+c}{2}\id&0
\\y^t&-u^t&0&-B^t+\frac{a+c}{2}\id\end{pmatrix},
\end{equation*}
where $\id$ denotes the $n\x n$ identity matrix. 

\begin{prop*}
  The map $i:Q\to\tilde P$ is an injective group homomorphism and the
  pair $(i,\al)$ satisfies conditions (1)--(3) of Proposition
  \ref{3.1}. Hence it gives rise to an extension functor from Cartan
  geometries of type $(G,Q)$ to Cartan geometries of type $(\tilde
  G,\tilde P)$.
\end{prop*}
\begin{proof}
All these facts are verified by straightforward computations, some of
which are a little tedious. 
\end{proof}

\subsection{Regularity and normality}\label{3.6}
We next have to discuss the conditions on the curvature of a Cartan
connection which were used in Theorems \ref{2.1} and \ref{2.3}.  If
$G$ is a semisimple group and $P\subset G$ is parabolic, then one can
identify $(\frak g/\frak p)^*$ with $\frak p_+$, the sum of all
positive grading components, via the Killing form, see \cite[Lemma 3.1]{Yam}.
Hence we can view the curvature function defined in \ref{3.3} as
having values in $\La^2\frak p_+\otimes\frak g$. Via the gradings of
$\frak p_+$ and $\frak g$, this space is naturally graded, and the
Cartan connection $\om$ is called \textit{regular} if its curvature
function has values in the part of positive homogeneity. Otherwise
put, if $X\in\frak g_i$ and $Y\in\frak g_j$, then $\ka(u)(X+\frak
p,Y+\frak p)\in \frak g_{i+j+1}\oplus\dots\oplus\frak g_k$.

Recall that a Cartan geometry is torsion free, if and only if $\ka$
has values in $\La^2\frak p_+\otimes\frak p$. Since elements of $\frak
p_+$ have strictly positive homogeneity, this subspace is contained in
the part of positive homogeneity, and any torsion free Cartan geometry
is automatically regular. Hence regularity should be viewed as a
condition which avoids particularly bad types of torsion.

On the other hand, there is a natural map $\partial^*:\La^2\frak
p_+\otimes\frak g\to \frak p_+\otimes\frak g$ defined by
$$
\partial^*(Z\wedge W\otimes A):=-W\otimes [Z,A]+Z\otimes
[W,A]-[Z,W]\otimes A 
$$ 
for decomposable elements. This is the differential in the standard
complex computing the Lie algebra homology of $\frak p_+$ with
coefficients in the module $\frak g$. This map is evidently
equivariant for the natural $P$--action, so in particular,
$\ker(\partial^*)\subset \La^2\frak p_+\otimes\frak g$ is a
$P$--submodule. The Cartan connection $\om$ is called \textit{normal}
if and only if its curvature has values in this submodule.

To proceed with the program set out in the end of \ref{3.4} we next
have to analyze the map $\Ps_\al:\La^2(\tg/\tp)\to\tg$ introduced in
\ref{3.3} in the special case of the pair $(i,\al)$ from \ref{3.5}. As
a linear space, we may identify $\tg/\tp$ with
$\tg_-=\tg_{-1}^E\oplus\tg_{-1}^V\oplus\tg_{-2}$. Note that using
brackets in $\tg$, we may identify $\tg_{-1}^V$ with
$\tg_1^E\otimes\tg_{-2}$ if necessary. We will view $\tg_{-2}$ as
$\Bbb R^{2n}=\Bbb R^n\oplus\Bbb R^n$ and correspondingly write
$X\in\tg_{-2}$ as $(X_1,X_2)$. By $\langle\ ,\ \rangle$ we denote the
standard inner product on $\Bbb R^n$.

\begin{lemma*}
Viewing $\Ps_\al$ as an element of $\La^2(\tg_-)^*\otimes\tg$, it lies in
the subspace $(\tg_{-1}^V)^*\wedge(\tg_{-2})^*\otimes\tg_0$. Denoting by
$W_0\in\tg_1^E$ the element whose unique nonzero entry is equal to $1$,
the trilinear map $\tg_{-2}\x \tg_{-2}\x \tg_{-2}\to\tg_{-2}$
defined by $(X,Y,Z)\mapsto [\Ps_{\al}(X,[Y,W_0]),Z]$ is (up to a nonzero
multiple) the complete symmetrization of the map $(X,Y,Z)\mapsto
\langle X_1,Y_2\rangle\begin{pmatrix} Z_1\\ -Z_2\end{pmatrix}$.      
\end{lemma*}
\begin{proof}
Let $x\in\tg_{-1}^E$ be the element whose unique nonzero entry is equal
to $1$. Then an arbitrary element of $\tg_-$ can be written uniquely as
$X+[Y,W_0]+ax$ for $X,Y\in\tg_{-2}$ and $a\in\Bbb R$. From the
definition of $\al$ in \ref{3.5} we obtain
$$
\al\begin{pmatrix} 0 & -Y_2^t & 0\\ X_1 & 0 & Y_1\\ a & X_2^t &
    0\end{pmatrix}=
\begin{pmatrix}
  0 & 0 & -\tfrac12 Y_2^t & \tfrac12 Y_1^t\\
a & 0 & \tfrac12 X_2^t & -\tfrac12 X_1^t \\
X_1 & Y_1 &0 &0\\ X_2 & Y_2 & 0 &0
\end{pmatrix}, 
$$
so this is congruent to $X+[Y,W_0]+ax$ modulo $\tp$. Using this,
one can now insert into the defining formula for $\Ps_\al$ from
\ref{3.3} and compute directly that the result always has values in
$\tg_0$, and indeed only in the lower right $2n\x
2n$ block. Moreover, all the entries in that block are made up from
bilinear expressions involving one entry from $\tg_{-2}$ and one entry
from $\tg_{-1}^V$, so we see that $\Ps_\al\in
(\tg_{-2})^*\wedge(\tg_{-1}^V)^*\otimes\tg_0$. 

For $X,Y\in\tg_{-2}$, one next computes that the only nonzero block in
$\Ps_{\al}(X,[Y,W_0])$ (which is a $2n\x 2n$--matrix) is explicitly given by
$$
\tfrac{1}{2}\begin{pmatrix}
X_1Y_2^t+Y_1X_2^t+(Y_2^tX_1+X_2^tY_1)\id & X_1Y_1^t+Y_1X_1^t\\
-X_2Y_2^t-Y_2X_2^t & -Y_2X_1^t-X_2Y_1^t-(Y_2^tX_1+X_2^tY_1)\id
\end{pmatrix}.
$$ 
To obtain $[\Ps_{\al}(X,[Y,W_0]),Z]\in\tg_{-2}$ for another element
$Z\in\tg_{-2}$, we now simply have to apply this matrix to
$\begin{pmatrix}Z_1\\ Z_2\end{pmatrix}$. Taking into account that
$\langle v,w\rangle=v^tw=w^tv$ for $v,w\in\Bbb R^n$ we obtain half the
sum of all cyclic permutations of
$$
(\langle X_1,Y_2\rangle+\langle Y_1,X_2\rangle)
\begin{pmatrix} Z_1\\ -Z_2\end{pmatrix},
$$
which is three times the total symmetrization of
$(X,Y,Z)\mapsto\langle X_1,Y_2\rangle\begin{pmatrix} Z_1\\
  -Z_2\end{pmatrix}$. 
\end{proof}

Using this we can now complete the first part of the program outlined
in the end of \ref{3.4}:
\begin{thm*}
The extension functor associated to the pair $(i,\al)$ from \ref{3.5}
maps locally flat Cartan geometries of type $(G,Q)$ to torsion free
(and hence regular), normal parabolic geometries of type $(\tilde
G,\tilde P)$. 
\end{thm*}
\begin{proof}
 Let $(p:\Cal G\to N,\om)$ be a locally flat Cartan geometry of type
 $(G,Q)$. This means that $\om$ has trivial curvature, so by
 Proposition \ref{3.3}, the curvature function $\tilde\ka$ of the
 parabolic geometry $(\Cal G\x_i\tilde P,\tilde\om_\al)$ has the
 property that 
$$
\tilde\ka(j(u))=\Ps_{\al}:\La^2(\tg/\tp)\to\tg,
$$
 where $j:\Cal G\to\Cal G\x_i\tilde P$ is the natural map. By the
 lemma above, $\tilde\ka(j(u))$ has values in $\tg_0\subset\tp$, and
 since having values in $\tp$ is a $\tilde P$--invariant property,
 torsion freeness follows.
 
 Similarly, since $\ker(\partial^*)$ is a $\tilde P$--submodule in
 $\La^2(\tg/\tp)^*\otimes\tg$, it suffices to show that
 $\partial^*(\Ps_{\al})=0$ to complete the proof of the theorem. This
 may be checked by a direct computation, but there is a more
 conceptual argument: Tracefree matrices in the lower right $2n\x 2n$
 block of $\tg_0$ form a Lie subalgebra isomorphic to
 $\frak{sl}(2n,\Bbb R)$ which acts on each of the spaces
 $\La^k(\tg/\tp)^*\otimes\tg$. Hence we may decompose each of them
 into a direct sum of irreducible representations. Since $\partial^*$
 is a $\tilde P$--homomorphism, it is equivariant for this action of
 $\frak{sl}(2n,\Bbb R)$, and hence it can be nonzero only between
 isomorphic irreducible components.
 
 In the proof of the lemma we have noted that $\tg_{-2}$ is the
 standard representation of $\frak{sl}(2n,\Bbb R)$, so the explicit
 formula for $\Ps_\al$ shows that it sits in a component isomorphic to
 $S^3\Bbb R^{2n*}\otimes\Bbb R^{2n}$. There is a unique trace from
 this representation to $S^2\Bbb R^{2n*}$, and the kernel of this is
 well known to be irreducible. One immediately checks that
 $(\tg/\tp)^*\otimes\tg$ cannot contain an irreducible component
 isomorphic to the kernel of this trace. Hence we can finish the proof
 by showing that $\Ps_\al$ lies in the kernel of that trace, which is
 a simple direct computation.  
\end{proof}

This has a nice immediate application:
\begin{cor*}
  Consider the homogeneous model $G\to G/P$ of Lagrangean contact
  structures. Then the resulting path geometry of chains is non--flat
  and hence not locally isomorphic to $\tilde G/\tilde P$, but its
  automorphism group contains $G$. In particular, for each $n\geq 1$,
  we obtain an example of a non--flat torsion free path geometry on a
  manifold of dimension $2n+1$ whose automorphism group has dimension
  at least $n^2+4n+3$.
\end{cor*}

\subsection*{Remark}
  (1) In \cite{Grossman}, the author directly constructed a torsion
  free path geometry from the homogeneous model of
  three--dimensional Lagrangean contact structures. This construction
  was one of the motivations for this paper and one of the
  guidelines for the right choice of the pair $(i,\al)$. The other
  main guideline for this choice are the computations needed to
  show that $\Ps_\al$ has values in $\tg_0$. 

\noindent
(2) We shall see later that in the situation of the corollary, the
dimension of the automorphism group actually equals the dimension of
$G$. In particular, for $n=1$, one obtains a non--flat path geometry
on a three manifold with automorphism group of dimension 8. To our
knowledge, this is the maximal possible dimension for the automorphism
group of a non--flat path geometry in this dimension.

Via the interpretation of path geometries in terms of systems of
second order ODE's, we obtain examples of nontrivial systems of such
ODE's with large automorphism groups.  

\subsection{More on curvatures of regular normal geometries}\label{3.7}
We have completed half of the program outlined in the end of \ref{3.4}
at this point: Theorem \ref{3.6} shows that the extension functor
associated to the pair $(i,\al)$ defined in \ref{3.5} produces the
regular normal parabolic geometry determined by the path geometry of
chains for locally flat Lagrangean contact structures. In view of
Theorem \ref{3.4} and Lemma \ref{3.1} this pins down the pair
$(i,\al)$ up to equivalence and hence the associated extension functor
up to isomorphism. 

Hence it only remains to clarify under which conditions on a
Lagrangean contact structure this extension procedure produces a
regular normal parabolic geometry. This then tells us the most general
situation in which a direct relation (as discussed in \ref{2.4} and
\ref{3.1}) between the two parabolic geometries can exist. As it can
already be expected from the case of the Fefferman construction (see
\cite{Cap:Fefferman}) this is a rather subtle question. Moreover, the
result cannot be obtained by algebraically comparing the two
normalization conditions, but one needs more information on the
curvature of regular normal and torsion free normal geometries. In
particular, the proof of part (2) of the Lemma below needs quite a lot
of deep machinery for parabolic geometries.

As discussed in \ref{3.6}, the curvature function of a parabolic
geometry of type $(G,P)$ has values in $\La^2\frak p_+\otimes\frak
g$. Since both $\frak p_+$ and $\frak g$ are graded, there is a
natural notion of homogeneity on this space. While being of some fixed
homogeneity is not a $P$--invariant property, the fact that all
nonzero homogeneous components have at least some given homogeneity is
$P$--invariant. This is used in the definition of regularity in
\ref{3.6}, which simply says that all nonzero homogeneous components
are in positive homogeneity. 

The map $\partial^*$ used in the definition of normality in \ref{3.6}
actually extends to a family of maps $\partial^*:\La^\ell\frak
p_+\otimes\frak g\to\La^{\ell-1}\frak p_+\otimes\frak g$. These are
the differentials in the standard complex computing the Lie algebra
homology $H_*(\frak p_+,\frak g)$. By definition, the curvature
function $\ka$ of a normal parabolic geometry of type $(G,P)$ has
values in $\ker(\partial^*)\subset\La^2\frak p_+\otimes\frak g$. Hence
we can naturally project to the quotient to obtain a function $\ka_H$
with values in $\ker(\partial^*)/\im(\partial^*)=H_2(\frak p_+,\frak
g)$. Equivariancy of $\ka$ implies that $\ka_H$ can be viewed as a
smooth section of the bundle $\Cal G\x_PH_2(\frak p_+,\frak g)$. This
section is called the \textit{harmonic curvature} of the normal
parabolic geometry. It turns out (see \cite{CSS-BGG}) that $P_+$ acts
trivially on $H_*(\frak p_+,\frak g)$, so this bundle admits a direct
interpretation in terms of the underlying structure. As we shall see
below, this bundle is algorithmically computable.

Now from \ref{3.6} we know that $\frak p_+\cong(\frak g/\frak p)^*$ as
a $P$--module, and since $\frak g_-\subset\frak g$ is a complementary
subspace (and $G_0$--module) to $\frak p\subset\frak g$ we can
identify $\frak p_+$ with $(\frak g_-)^*$ as a $G_0$--module. Hence we
can also view the spaces $\La^\ell\frak p_+\otimes\frak g$ as
$L(\La^\ell\frak g_-,\frak g)$, which are the chain spaces in the
standard complex computing the Lie algebra cohomology of $\frak g_-$
with coefficients in $\frak g$. The differentials
$\partial:L(\La^\ell\frak g_-,\frak g)\to L(\La^{\ell+1}\frak
g_-,\frak g)$ in that complex turn out to be adjoint to the maps
$\partial^*$ with respect to a certain inner product.

Hence we obtain an algebraic Hodge theory on each of the spaces
$\La^\ell\frak p_+\otimes\frak g$, with algebraic Laplacian
$\square=\partial^*\o\partial+\partial\o\partial^*$. This construction
is originally due to Kostant (see \cite{Kostant}), whence $\square$ is
usually called the Kostant Laplacian. The kernel of $\square$ is a
$G_0$--submodule called the \textit{harmonic subspace} of
$\La^\ell\frak p_+\otimes\frak g$. Kostant's version of the
Bott--Borel--Weil theorem in \cite{Kostant} gives a complete
algorithmic description of the $G_0$--module $\ker(\square)$. By the
Hodge decomposition, $\ker(\square)$ is isomorphic to the homology
group of the appropriate dimension.

We will need two general facts about the curvature of regular normal
respectively torsion free normal parabolic geometries in the sequel: 
\begin{lemma*}
Let $(p:\Cal G\to M,\om)$ be a regular normal parabolic geometry of
type $(G,P)$ with curvature functions $\ka:\Cal G\to \La^2\frak
p_+\otimes\frak g$ and $\ka_H:\Cal G\to H_2(\frak p_+,\frak g)$. 
Then we have:

\noindent
(1) The lowest nonzero homogeneous component of $\ka$ has values in
    the subset $\ker(\square)\subset\La^2\frak p_+\otimes\frak g$. 

\noindent
(2) Suppose that $(p:\Cal G\to M,\om)$ is torsion free and that
    $E_0\subset\ker(\square)\subset\La^2\frak p_+\otimes\frak g$ is a
    $G_0$--submodule such that $\ka_H$ has values in the image of $E_0$
    under the natural isomorphism $\ker(\square)\to H_2(\frak
    p_+,\frak g)$ (induced by projecting
    $\ker(\square)\subset\ker(\partial^*)$ to the quotient). Then $\ka$
    has values in the $P$--submodule of $\La^2\frak p_+\otimes\frak g$
    generated by $E_0$. 
\end{lemma*}
\begin{proof}
(1) is an application of the Bianchi identity, which goes back to
    \cite{Tan}, see also \cite[Corollary 4.10]{C-S}. (2) is
    proved in \cite[Corollary 3.2]{C-tw}.     
\end{proof}

The final bit of information we need is the explicit form of
$\ker(\square)$ for the pairs $(\frak g,\frak p)$ and $(\tg,\tp)$
corresponding to Lagrangean contact structures on manifolds of
dimension $2n+1$ respectively path geometries in dimension $4n+1$.
Obtaining the explicit description of the irreducible components of
these submodules is an exercise in the application of Kostant's
results from \cite{Kostant} and the algorithms from the book
\cite{BE}, see also \cite{C-tw}. The results are listed in the tables
below. The first column contains the homogeneity of the component and
the second column contains the subspace that it is contained in. The
actual component is always the highest weight part in that subspace,
so in particular, it lies in the kernel of all traces one can form.

\begin{tabular}{cc}
\parbox[c][][c]{0.45\textwidth}{
\begin{center}
$(\frak g,\frak p)$, $n=1$ \\[6pt] 
\begin{tabular}{|c|c|}
  \hline homog. & \parbox[c][1.4\totalheight][c]{55pt}{contained in} \\
  \hline 4 & \parbox[c][1.4\totalheight][c]{55pt}{$\g_1^R\wedge\g_2\otimes\g_1^R$}\\
  \hline 4 & \parbox[c][1.4\totalheight][c]{55pt}{$\g_1^L\wedge\g_2\otimes\g_1^L$}\\
  \hline
\end{tabular}
\end{center}}&
\parbox[c][][c]{0.45\textwidth}{
\begin{center}
$(\frak g,\frak p)$, $n>1$ \\[6pt]
\begin{tabular}{|c|c|}
  \hline homog. & \parbox[c][1.4\totalheight][c]{55pt}{contained in}\\
  \hline 2 & \parbox[c][1.4\totalheight][c]{55pt}{$\g_1^L\wedge\g_1^R\otimes\g_0$} \\
  \hline 1 & \parbox[c][1.4\totalheight][c]{55pt}{$\La^2\g_1^L\otimes\g_{-1}^R$} \\
  \hline 1 & \parbox[c][1.4\totalheight][c]{55pt}{$\La^2\g_1^R\otimes\g_{-1}^L$} \\
  \hline
\end{tabular}
\end{center}}
\end{tabular}

\bigskip

\begin{tabular}{cc}
\parbox[c][][c]{0.45\textwidth}{
\begin{center}
$(\tg,\tp)$, $n=1$ \\[6pt] 
\begin{tabular}{|c|c|c|}
  \hline homog. & \parbox[c][1.4\totalheight][c]{55pt}{contained in}\\
  \hline 3 & \parbox[c][1.4\totalheight][c]{55pt}{$\tg_1^V\wedge\tg_2\otimes\tg_0$} \\
  \hline 2 & \parbox[c][1.4\totalheight][c]{55pt}{$\tg_1^E\wedge\tg_2\otimes\tg_{-1}^V$} \\
  \hline 1 & \parbox[c][1.4\totalheight][c]{55pt}{$\La^2\tilde\g_1^V\otimes\tilde\g_{-1}^E$} \\
  \hline
\end{tabular}
\end{center}}&
\parbox[c][][c]{0.45\textwidth}{
\begin{center}
$(\tg,\tp)$, $n>1$ \\[6pt]
\begin{tabular}{|c|c|c|}
  \hline homog. & \parbox[c][1.4\totalheight][c]{55pt}{contained in} \\
  \hline 3 & \parbox[c][1.4\totalheight][c]{55pt}{$\tg_1^V\wedge\tg_2\otimes\tg_0$} \\
  \hline 2 & \parbox[c][1.4\totalheight][c]{55pt}{$\tg_1^E\wedge\tg_2\otimes\tg_{-1}^V$} \\
  \hline 0 & \parbox[c][1.4\totalheight][c]{55pt}{$\La^2\tg_1^V\otimes\tg_{-2}$} \\
  \hline
\end{tabular}
\end{center}}
\end{tabular}

\subsection{}\label{3.8}
We are now ready to prove the main result of this article:
\begin{thm*}
Let $(p:\Cal G\to M,\om)$ be a regular normal parabolic geometry of
type $(G,P)$ and let $(\tcg:=\Cal G\x_Q\tilde P\to\Cal
P_0(TM),\tilde\om_\al)$ be the parabolic geometry obtained using the
extension functor associated to the pair $(i,\al)$ defined in
\ref{3.5}. Then this geometry is regular and normal if and only if
$(p:\Cal G\to M,\om)$ is torsion free. 
\end{thm*}
\begin{proof}
  We first prove necessity of torsion freeness. From the tables in
  \ref{3.7} we see that for $n=1$ a regular normal parabolic geometry
  of type $(G,P)$ is automatically torsion free, so we only have to
  consider the case $n>1$. If $\tilde\om_\al$ is regular and normal,
  then all nonzero homogeneous components of $\tilde\ka$ are
  homogeneous of positive degrees. The table in \ref{3.7} shows that
  then the homogeneity is at least two, and by part (1) of Lemma
  \ref{3.7} the homogeneous component of degree two sits in the
  subspace $\tg_1^E\wedge\tg_2\otimes \tg_{-1}^V$. In particular, 
  for any $\tilde u\in\tilde\G$, the
  restriction of $\tilde\ka(\tilde u)$ to $\La^2\tg_{-2}$ is homogeneous of
  degree at least three, which implies that $\tilde\ka(\tilde u)$ has values
  in $\tg_{-1}\oplus\tp$, i.e.~for the natural projection $\pi:\tg\to
  \tg/(\tg_{-1}\oplus\tp)$ we get $\pi\o\tilde\ka(\tilde u)=0$.
  
  Using the notation of the proof of Lemma \ref{3.6}, consider two
  elements $X,Y\in\tg_{-2}$. From that proof, we see that
$$
(\pi\o\tilde\ka(j(u)))(X,Y)=(\pi\o\al\o\ka(u))
  \left(\begin{pmatrix} 0 & 0 & 0\\ X_1 & 0 & 0\\
  0 & X_2^t & 0\end{pmatrix},\begin{pmatrix} 0 & 0 & 0\\ Y_1 & 0 & 0\\
  0 & Y_2^t & 0\end{pmatrix}\right).
$$
By regularity,
$\ka(u)(\La^2\fg_{-1})\subset\fg_{-1}\oplus\frak p$. From the
definition in \ref{3.5} it is evident that $\al$ induces a linear
isomorphism $\fg/(\fg_{-2}\oplus\frak p)\to
\tg/(\tg_{-1}\oplus\tp)$. Hence we conclude that if
$\tilde\om_\al$ is regular and normal, then
$\ka(u)(\La^2\fg_{-1})\subset\frak p$. From the table in \ref{3.7} we
see that this implies that the homogeneous component of degree one of
$\ka$ has to vanish identically, and then further that the homogeneous
component of degree two has values in $\frak p$. Since $\La^2\frak
g_{-2}=0$, components of homogeneity at least three automatically have
values in $\frak p$, so we see that $\om$ is torsion free.

To prove sufficiency, we first need two facts on the curvature function
$\ka$ of a torsion free normal parabolic geometry of type $(G,P)$. On
the one hand, the map $\partial^*$ as defined in \ref{3.6} can be
written as the sum $\partial_1^*+\partial_2^*$ of two $P$--equivariant
maps, with $\partial^*_1$ corresponding to the first two summands and
$\partial_2^*$ corresponding to the last summand in the definition. We
claim that $\ka$ has values in the kernels of both operators
$\partial^*_i$. On the other hand, one easily verifies that the
subspace $\widehat{\frak p}\subset\frak p$ formed by all matrices of
the form $\begin{pmatrix}0 & u & d\\ 0 & B & v\\ 0 & 0 & 0
\end{pmatrix}$ is a $P$--submodule. (Indeed, this is the preimage in
$\frak p$ of the semisimple part of the reductive algebra $\frak
g_0=\frak p/\frak p_+$.) Our second claim is that
$\ka(u)(X,Y)\in\widehat{\frak p}$ for all $u\in\Cal G$ and all 
$X,Y$.

To prove both claims, it suffices to show that $\ka$ has values in the
$P$--submodule $\La^2_0\frak p_+\otimes\widehat{\frak
  p}\subset\La^2\frak p_+\otimes\frak p$. Here $\La^2_0\frak p_+$ is
the kernel of the $P$--homomorphism $\La^2\frak p_+\to\frak p_+$
defined by the Lie bracket on $\frak p_+$, so $ \La^2_0\frak
p_+\otimes\frak g=\ker(\partial^*_2)$. 

In the case $n=1$, this is evident, since from the table in \ref{3.7}
we see that the lowest nonzero homogeneous component of $\ka(u)$ is
of degree $4$, vanishes on $\La^2\frak g_{-1}$ and has values in
$\frak p_+$. For homogeneous components of higher degree, these two
properties are automatically satisfied, and we conclude that
$\ka(u)\in \fg_1\wedge\fg_2\otimes\frak p_+\subset\La^2_0\frak
p_+\otimes\widehat{\frak p}$.

In the case $n>1$, we see from the table in \ref{3.7} that by torsion
freeness the lowest homogeneous component of $\ka(u)$ must be of
homogeneity $2$. By part (1) of Lemma \ref{3.7} it has values in
$\ker(\square)\subset\La^2\fg_1\otimes\fg_0$. Since this component
of $\ker(\square)$ is a highest weight part, it lies in the kernel of
all possible traces, and hence it must be contained in the tensor
product of $\La^2\fg_1\cap \La^2_0\frak p_+$ with the semisimple part
of $\frak g_0$. Hence $\ker(\square)$ is contained in the
$P$--submodule $\La^2_0\frak p_+\otimes\widehat{\frak p}$ so, by part
(2) of Lemma \ref{3.7}, the curvature function $\ka$ has values in
that submodule. 

In view of Proposition \ref{3.3} and the proof of Theorem \ref{3.6},
to prove that $\tilde\om_\al$ is regular and normal, it suffices to
verify that the map $F(u):\La^2\tg_-\to\tg$ defined by
$F(u)(X,Y):=\al(\ka(u)(\underline{\al}^{-1}(X),\underline{\al}^{-1}(Y)))$
lies in the kernel of $\partial^*$ for all $u\in\Cal G$. To compute
$\partial^* F(u)$, it is better to view $F(u)$ as an element of
$\La^2\tp_+\otimes\tg$, and we want to relate this to $\ka(u)$, viewed
as an element of $\La^2\frak p_+\otimes\frak g$. Therefore, we have to
compute the map $\ph:\frak p_+\to\tp_+$, which is dual to the
composition of the canonical projection $\frak g/\frak q\to\frak
g/\frak p$ with $\underline{\al}^{-1}:\tg/\tp\to\frak g/\frak q$,
since by construction $F(u)=(\La^2\ph\otimes\al)(\ka(u))$.  Recall
that the duality between $\frak g/\frak p$ and $\frak p_+$ (and
likewise for the other algebra) is induced by the Killing form. Since
the Killing form of a simple Lie algebra is uniquely determined up to
a nonzero multiple by invariance, we may as well use the trace form on
both sides, which leads to a nonzero multiple of $\ph$. But then the
computation is very easy, showing that
$$
\ph\begin{pmatrix} 0 & Z & \ps \\ 0 & 0 & W\\
    0&0&0\end{pmatrix}= 
\begin{pmatrix}0 & \ps & Z & W^t\\ 0&0&0&0\\ 0&0&0&0\\  0&0&0&0\end{pmatrix}.
$$
In particular, $\ph(\frak p_+)\subset \tg_1^E\oplus\tg_2$, which
implies that $\partial^*_2(F(u))=0$ for all $u$.

On the other hand, the formula for $\al$ from \ref{3.5} shows that
$\al(\widehat{\frak p})\subset
\tg_{-1}^V\oplus\tg_0\oplus\tg_1^E\oplus\tg_2$, and the
$\tg_0$--component is contained in the bottom right $2n\x 2n$ block.
This shows that for $Z\in\frak p_+$ and $A\in\widehat{\frak p}$ we
have $[\ph(Z),\al(A)]\in \tg_1^E\oplus\tg_2$. One immediately verifies
directly that the $\tg_1^E$--component of $[\ph(Z),\al(A)]$ equals the
$\tg_1^E$--component of $\al([Z,A])$, while the $\tg_2$--component of
$[\ph(Z),\al(A)]$ equals twice the $\tg_2$--component of $\al([Z,A])$.
From the definition of $\partial^*_1$ we now conclude that
$\La^2\ph\otimes\al$ maps $\ker(\partial_1^*)$ to
$\ker(\partial_1^*)$, so we also get $\partial^*_1(F(u))=0$ for all
$u$.
\end{proof}

\section{Applications}\label{4}
For torsion free Lagrangean contact structures, Theorem \ref{3.8}
provides us with an explicit description of the parabolic geometry
determined by the path geometry of chains. In particular, we obtain an
explicit formula for the Cartan curvature which is the basis for the
applications discussed in this section. The main result is that one
can essentially reconstruct the torsion free Lagrangean contact
structure from the harmonic curvature of this parabolic geometry. In
particular, this implies that a contact diffeomorphism which maps
chains to chains has to either preserve or swap the subbundles
defining the Lagrangean contact structure. On the way, we can prove
that chains can never be described by linear connections and that only
locally flat Lagrangean contact structures give rise to torsion free
path geometries of chains.

\subsection{Decomposing the Cartan curvature}\label{4.1}
For a torsion free Lagrangean contact structure with curvature $\ka$,
the curvature $\tilde\ka$ of the normal Cartan connection associated
to the path geometry of chains is determined by the formula from
Proposition \ref{3.3}, which holds on $j(\Cal G)\subset\Cal
G\x_i\tilde P$. In this formula, there are two terms, one of which
depends on $\ka$ while the other one only comes from the map $\al$.
Our main task is to extract parts of $\tilde\ka$ which only depend on
one of the two terms. The difficulty is that this has to be done in a
geometric way without knowing the subset $j(\Cal G)$ in advance.

The curvature function $\tilde\ka$ has values in the $P$--module
$\La^2\tp_+\otimes\tg$, and using the map $\ph$ from the proof of
Theorem \ref{3.8}, the formula from Proposition \ref{3.3} reads as
$\tilde\ka(j(u))=(\La^2\ph\otimes\al)(\ka(u))+\Ps_\al$. Now $\tp_+$
contains the $P$--invariant subspace $\tg_2$. Correspondingly, we
obtain $P$--invariant subspaces
$\La^2\tg_2\subset\tp_+\wedge\tg_2\subset\La^2\tp_+$. In the proof of
Theorem \ref{3.8}, we have seen that $\ph$ has values in
$\tg_1^E\oplus\tg_2$, whence $\La^2\ph$ has values in
$\tp_+\wedge\tg_2$. From Lemma \ref{3.6} we know that
$\Ps_\al\in\tp_+\wedge\tg_2\otimes\tg$, so we conclude that
$\tilde\ka(j(u))$ lies in this $\tilde P$--submodule. By equivariancy,
all values of the curvature function lie in
$\tp_+\wedge\tg_2\otimes\tg\subset\La^2\tp_+\otimes\tg$.

On the quotient $\tp_+/\tg_2$, the subgroup $\tilde P_+\subset\tilde
P$ acts trivially, so we can identify this quotient with the $\tilde
G_0$--module $\tg_1=\tg_1^E\oplus\tg_1^V$. Correspondingly, we get
$\tilde P$--equivariant projections 
\begin{gather*}
\pi^E:\tp_+\wedge\tg_2\otimes\tg\to \tg_1^E\wedge\tg_2\otimes\tg \\
\pi^V:\tp_+\wedge\tg_2\otimes\tg\to \tg_1^V\wedge\tg_2\otimes\tg.
\end{gather*}
From the description of the image of $\ph$ in the proof of Theorem
\ref{3.8} we conclude that
$(\La^2\ph\otimes\al)(\ka(u))\in\ker(\pi^V)$. On the other hand, Lemma
\ref{3.6} in particular shows that $\pi^V(\Ps_\al)\neq 0$ and
$\Ps_\al\in\ker(\pi^E)$.

\begin{thm*}
Let $(M,L,R)$ be a torsion free Lagrangean contact structure.

\noindent
(1) There is no linear connection on the tangent bundle $TM$ which has
    the chains among its geodesics. 

\noindent
(2) The parabolic geometry associated to the path geometry of chains
    on $\tilde M=\Cal P_0(TM)$ is torsion free if and only if
    $(M,L,R)$ is locally flat, i.e.~locally isomorphic to the
    homogeneous model $G/P$. 
\end{thm*}
\begin{proof}
  (1) Suppose that $\nabla$ is a linear connection on $TM$ whose
  geodesics in directions transverse to $L\oplus R$ are
  parametrizations of the chains. Since symmetrizing a connection does
  not change the geodesics, we may without loss of generality assume
  that $\nabla$ is torsion free. Then we can look at the associated
  projective structure $[\nabla]$ on $M$ and use the machinery of
  correspondence space from \cite{C-tw}. The fact that the geodesics
  of $\nabla$ are the chains exactly means that the path geometry of
  chains on $\tilde M$ is isomorphic to an open subgeometry of the
  correspondence space $\Cal C(M,[\nabla])$, see 4.7 of \cite{C-tw}.
  In particular, the Cartan curvature $\tilde\ka$ is the restriction
  of the curvature of this correspondence space. By \cite[Proposition
  2.4]{C-tw} this curvature has the property that it vanishes upon
  insertion of one tangent vector contained in the vertical bundle of
  $\tilde M\to M$. But this contradicts the fact that
  $\pi^V\o\tilde\ka\neq 0$ we have observed above.

\noindent
(2) By Theorem \ref{3.6}, the path geometry of chains associated to a
locally flat Lagrangean contact structure is torsion free. Conversely,
if the Cartan connection $\tilde\om$ is torsion free, then according
to part (1) of Lemma \ref{3.7} and the tables in \ref{3.7}, the lowest
nonzero homogeneous component of $\tilde\ka$ must be of degree at
least three, and the harmonic curvature must have values in
$\tg_1^V\wedge\tg_2\otimes\tg_0\subset\ker(\pi^E)$.  By part (2) of
Lemma \ref{3.7} the whole curvature $\tilde\ka$ has values in
$\ker(\pi^E)$. Above, we have observed that $\Ps_\al\in\ker(\pi^E)$ so
we conclude that for each $u\in\Cal G$ we get
$\pi^E\o (\La^2\ph\otimes\al)(\ka(u))=0$. 

In the proof of Theorem \ref{3.8} we see that $\ph$ is a linear
isomorphism $\frak p_+\to \tg_1^E\oplus\tg_2$, and hence
$\tg_1^E\wedge\tg_2$ is contained in the image of $\La^2\ph$. Hence we
conclude that $\al\o\ka(u)=0$ and since $\al$ is injective, the result
follows.
\end{proof}

\subsection{Harmonic curvature}\label{4.2}
We have discussed the definition of harmonic curvature already in
\ref{3.7}. Let $\pi_H$ be the natural projection from
$\ker(\partial^*)\subset\La^2\tp_+\otimes\tg$ to the quotient
$\ker(\partial^*)/\im(\partial^*)$. Since this is a $\tilde
P$--equivariant map, the composition
$\tilde\ka_H=\pi_H\o\tilde\ka:\tcg\to
\ker(\partial^*)/\im(\partial^*)$ defines a smooth section of the
associated bundle $\tcg\x_{\tilde P}\ker(\partial^*)/\im(\partial^*)$,
which is the main geometric invariant of the parabolic geometry
associated to the path geometry of chains. 

From \ref{3.7} we also know that $\tilde P_+$ acts trivially on the
quotient $\ker(\partial^*)/\im(\partial^*)$ and we may identify it
with the $\tilde G_0$--module
$\ker(\square)\subset\La^2\tp_+\otimes\tg$. From the table in
\ref{3.7}, we see that this module contains two irreducible components
in positive homogeneity, which are the highest weight components of
the subrepresentations $\tg_1^E\wedge\tg_2\otimes\tg_{-1}^V$
respectively $\tg_1^V\wedge\tg_2\otimes\tg_0$. Correspondingly, we
obtain decompositions $\pi_H=\pi^E_H+\pi_H^V$ and
$\tilde\ka_H=\tilde\ka_H^E+\tilde\ka_H^V$. 

\begin{lemma*}
Let $\pi^E$ and $\pi^V$ be the projections on
$\tp_+\wedge\tg_2\otimes\tg$ defined in \ref{4.1}. Then the
restriction of $\pi_H^E$ (respectively $\pi_H^V$) to
$\ker(\partial^*)\cap(\tp_+\wedge\tg_2\otimes\tg)$ factorizes through
$\pi^E$ (respectively $\pi^V$).
\end{lemma*}
\begin{proof}
  By Kostant's version of the Bott--Borel--Weil theorem, see
  \cite{Kostant}, the $\tilde G_0$--irreducible components contained
  in $\ker(\square)$ occur with multiplicity one, even within
  $\La^*\tp_+\otimes\tg$. To obtain $\pi^E$ and $\pi^V$, we used the
  projection $\tp_+\wedge\tg_2\otimes\tg\to\tg_1\wedge\tg_2\otimes\tg$
  with kernel $\La^2\tg_2\otimes\tg$. By the multiplicity one result
  and the fact that both components of $\ker(\square)$ are contained
  in $\tg_1\wedge\tg_2\otimes\tg$, there is no nonzero $\tilde
  G_0$--equivariant map $\La^2\tg_2\otimes\tg\to
  \ker(\partial^*)/\im(\partial^*)$. Hence each of the projections
  $\pi_H$, $\pi_H^E$ and $\pi_H^V$ factorizes through
  $\tg_1\wedge\tg_2\otimes\tg$. Looking at the resulting map for
  $\pi_H^E$, we see that again by multiplicity one, the subspace
  $\tg_1^V\wedge\tg_2\otimes\tg$ must be contained in the kernel, so
  we conclude that $\pi_H^E$ factorizes through $\pi^E$. In the same
  way one shows that $\pi_H^V$ factorizes through $\pi^V$.
\end{proof}

\begin{prop*}
Let $(M,L,R)$ be a torsion free Lagrangean contact structure, and let
$\tilde\ka_H=\tilde\ka_H^E+\tilde\ka_H^V$ be the harmonic curvature of
the regular normal parabolic geometry determined by the path geometry
of chains.

Then the function $\tcg\to \tg_1^V\wedge\tg_2\otimes\tg_0$
corresponding to $\tilde\ka_H^V$ is a nonzero multiple of the unique
equivariant extension of the constant function $\Ps_\al$ (compare with
Lemma \ref{3.6}) on $j(\Cal G)$.
\end{prop*}
\begin{proof}
We have to compute the function $\pi_H^V\o\tilde\ka$. By the lemma,
$\pi_H^V$ factorizes through the projection $\pi^V$ introduced in
\ref{4.1}, and from there we know that
$\pi_V(\tilde\ka(j(u)))=\pi_V(\Ps_\al)$. Hence we see that
$(\pi_H^V\o\tilde\ka)|_{j(\Cal G)}=\pi_H^V(\Ps_\al)$. Now $\Ps_\al\in
\tg_1^V\wedge\tg_2\otimes\tg_0$ by Lemma \ref{3.6}, and the values
even lie in the semisimple part of $\tg_0$, which may be identified
with $\frak{sl}(\tg_{-2})$. Evidently, $\tg_1^V\cong
\tg_{-1}^E\otimes\tg_2$ as a $\tilde G_0$--module, so we may interpret
$\Ps_\al$ as an element of
$\tg_{-1}^E\otimes(\otimes^3\tg_2)\otimes\tg_{-2}$. In Lemma \ref{3.6}
and the proof of Theorem \ref{3.6} we have seen that in this picture
$\Ps_\al$ lies in the irreducible component $ \tg_{-1}^E\otimes
(S^3\tg_2\otimes\tg_{-2})_0$, where the subscript denotes the
trace free part. Passing back to $ \tg_1^V\wedge\tg_2\otimes\tg_0$ this
exactly means that $\Ps_\al$ lies in the highest weight subspace,
which is the intersection with $\ker(\square)$. Now $\pi_H^V$
restricts to $\tilde G_0$--equivariant linear isomorphism on this
intersection, which implies the result.
\end{proof}

\noindent
\textbf{Remark.}  Similarly to the proof above, one shows that the
harmonic curvature component $\tilde\ka_H^E$ is the extension of a
component of $j(u)\mapsto(\La^2\ph\otimes\al)(\ka(u))$. Since we
explicitly know $\La^2\ph\otimes\al$, this can be used to obtain a
more explicit description of the second harmonic curvature
component. From part (2) of Theorem \ref{4.1} and \cite[4.7]{C-tw} we
see that vanishing of $\tilde\ka_H^E$ is equivalent to local flatness
of the original Lagrangean contact structure, so $\ka$ is completely
encoded in $\tilde\ka_H^E$.

\subsection{Passing to the underlying manifold}\label{4.3} 
The harmonic curvature component determined by the function
$\tilde\ka_H^V$ is a section of the bundle associated to
$\tg_1^V\wedge\tg_2\otimes\tg_0$. In the proof of Proposition
\ref{4.2} we have seen that we can replace that space by
$\tg_{-1}^E\otimes(\otimes^3\tg_2)\otimes\tg_2$. The corresponding
bundle is $E\otimes \otimes^3F^*\otimes F\to\tilde M$, where
$F:=T\tilde M/(E\oplus V)$. Since $E\subset TM$ is a line bundle, we
can view $\tilde\ka_H^V$ as a section of $\otimes^3F^*\otimes F$ which
is determined up to a nonzero multiple.

To relate this to the underlying manifold $M$, recall that $\tilde M$
is an open subset in the projectivized tangent bundle of $M$. A point
in $\tilde M$ is a line in some tangent space $T_xM$ that is
transversal to $L_x\oplus R_x$. We have noted in \ref{2.4} that
$TM\cong\Cal G\x_P\fg/\frak p$ and $T\tilde M\cong\Cal G\x_Q\fg/\frak
q$, and the tangent map of the projection $\pi:\tilde M\to M$
corresponds to the natural projection $\fg/\frak q\to\frak g/\frak
p$. Fix a point $\ell\in\pi^{-1}(x)$. Then for each $\xi\in T_xM$
there is a lift $\tilde\xi\in T_\ell\tilde M$ and we can consider the
class of $\tilde\xi$ in $F_\ell=T_\ell\tilde M/(E_\ell\oplus
V_\ell)$. Since $V_\ell$ is the vertical subbundle, this class is
independent of the choice of the lift and from the explicit
description of $T\pi$ we see that restricting to $L_x\oplus R_x$, we
obtain a linear isomorphism $L_x\oplus R_x\cong F_\ell$. 

Fixing $x$ and $\ell$ we therefore see that the harmonic curvature
component corresponding to $\tilde\ka_H^V$ gives rise to an element of
$\otimes^3(L_x\oplus R_x)^*\otimes (L_x\oplus R_x)$, which is
determined up to a nonzero multiple. To write down this map
explicitly, we first need the Levi bracket
$$
\Cal L:(L_x\oplus R_x)\x (L_x\oplus R_x)\to T_xM/(L_x\oplus R_x). 
$$ 
Since this has values in a one--dimensional space, we may view it as a
real valued bilinear map determined up to a nonzero multiple. Further,
we denote by $\Bbb J$ the almost product structure corresponding to
the decomposition $L\oplus R$. This means that $\Bbb J$ is the
endomorphism of $L\oplus R$ which is the identity on $L$ and minus the
identity on $R$. Using this we can now formulate:
\begin{lemma*}
  The element of $\otimes^3(L_x\otimes R_x)^*\otimes (L_x\oplus R_x)$
  obtained from $\tilde\ka_H^V$ above is (a nonzero multiple of) the
  complete symmetrization of the map
$$
(\xi,\eta,\zeta)\mapsto \Cal L(\xi,\Bbb J(\eta))\Bbb J(\zeta). 
$$
\end{lemma*}
\begin{proof}
This is a reinterpretation of the proof of Lemma \ref{3.6}. Observe
that $\Bbb J$ corresponds to the map $\begin{pmatrix} X_1\\
X_2\end{pmatrix}\mapsto\begin{pmatrix} X_1\\ -X_2\end{pmatrix}$ in the
notation there. Since $\Cal L$ corresponds to $[\ ,\
]:\fg_{-1}\x\fg_{-1}\to\fg_{-2}$, computing the bracket
$$
\left[\begin{pmatrix} 0 & 0 & 0\\ X_1 & 0 & 0 \\ 0 & X_2^t &
    0\end{pmatrix},\begin{pmatrix} 0 & 0 & 0\\ Y_1 & 0 & 0 \\
    0 & -Y_2^t & 0\end{pmatrix}\right],
$$
we see that the expression $\langle X_1,Y_2\rangle+\langle
Y_1,X_2\rangle$ in the proof of Lemma \ref{3.6} corresponds to
$\Cal L(\xi,\Bbb J(\eta))$. 
\end{proof}

\subsection{Reconstructing the Lagrangean contact
  structure}\label{4.4} 
Now we can finally show that the Cartan curvature of the path geometry
of chains can be used to (almost) reconstruct the Lagrangean contact
structure on $M$ that we have started from:

\begin{thm*}
  Let $(M,L,R)$ be a torsion free Lagrangean contact structure. Then
  for each $x\in M$, the subset $L_x\cup R_x\subset T_xM$ can be
  reconstructed from the harmonic curvature of the normal parabolic
  geometry associated to the path geometry of chains.
\end{thm*}
\begin{proof}
  In view of the results in \ref{4.2} and \ref{4.3} it suffices to
  show that $L_x\cup R_x$ can be recovered from the complete
  symmetrization $S$ of the map 
$$
(\xi,\eta,\zeta)\mapsto \Cal L(\xi,\Bbb J(\eta))\Bbb J(\zeta).
$$
First we see that $S(\xi,\xi,\xi)=0$ if and only if $\Cal L(\xi,\Bbb
J(\xi))=0$. Note that this is always satisfied for $\xi\in L_x\cup
R_x$. Fixing an element $\xi$ with this property, we see that
$$
S(\xi,\xi,\eta)=2\Cal L(\xi,\Bbb J(\eta))\Bbb J(\xi).
$$
By non--degeneracy of $\Cal L$, given a nonzero element $\xi$ we
can always find $\eta$ such that $\Cal L(\xi,\Bbb J(\eta))\neq
0$. Hence we see that $\xi$ is an eigenvector for $\Bbb J$ (which by
definition is equivalent to $\xi\in L_x\cup R_x$) if and only if
$S(\xi,\xi,\xi)=0$ and there is an element $\eta$ such that
$S(\xi,\xi,\eta)$ is a nonzero multiple of $\xi$.
\end{proof}

\begin{cor*}
  Let $(M,L,R)$ be a torsion free Lagrangean contact structure and let
  $f:M\to M$ be a contact diffeomorphism which maps chains to chains.
  Then either $f$ is an automorphism or an anti--automorphism of the
  Lagrangean contact structure. Here anti--automorphism means that
  $T_xf(L_x)=R_{f(x)}$ and $T_xf(R_x)=L_{f(x)}$ for all $x\in M$.
\end{cor*}
\begin{proof}
  By assumption, $f$ induces an automorphism $\tilde f$ of the path
  geometry of chains associated to $(M,L,R)$. This automorphism has to
  pull back the Cartan curvature $\tilde\ka$ and also the harmonic
  curvature $\ka_H$ to itself. From the theorem we conclude that this
  implies $T_xf(L_x\cup R_x)=L_{f(x)}\cup R_{f(x)}$, and this is only
  possible if $f$ is an automorphism or an anti--automorphism.
\end{proof}

\section{Partially integrable almost CR structures}\label{5}
What we have done for Lagrangean contact structures so far can be
easily adapted to deal with partially integrable almost CR
structure. We will only briefly sketch the necessary changes in this
section. 

\subsection{}\label{5.1}
A non--degenerate partially integrable almost CR structure on a smooth
manifold $M$ is given by a contact structure $H\subset TM$ together
with an almost complex structure $J$ on $H$ such that the Levi bracket
$\Cal L$ has the property that $\Cal L(J\xi,J\eta)=\Cal L(\xi,\eta)$
for all $\xi,\eta$. Then $\Cal L$ is the imaginary part of a
non--degenerate Hermitian form and we denote the signature of this
form by $(p,q)$. Such a structure of signature $(p,q)$ is equivalent
to a regular normal parabolic geometry of type $(G,P)$, where
$G=PSU(p+1,q+1)$ and $P\subset G$ is the stabilizer of a point in $\Bbb
CP^{n+1}$, $n=p+q$, corresponding to a null line, see
\cite[4.15]{C-S}. The group $G$ is the quotient of $SU(p+1,q+1)$ by
its center (which is isomorphic to $\Bbb Z_{n+2}$) and we will work
with representative matrices as before.

We will use the Hermitian form of signature $(p,q)$ on $\Bbb C^{n+1}$
corresponding to
$$
(z_0,\dots,z_{n+1})\mapsto z_0\bar z_{n+1}+z_{n+1}\bar z_0+
\textstyle\sum_{j=1}^p|z_j|^2-\textstyle\sum_{j=p+1}^n|z_j|^2.
$$
Then the decomposition on $\frak{sl}(n+2,\Bbb C)$ with block sizes 1,
$n$, and $1$ restricts to a contact grading on the Lie algebra $\frak
g$ of $G$. The explicit form for signature $(n,0)$ can be found in
\cite[4.15]{C-S}. In general, $\fg$ consists of all matrices of the
form 
$$
\begin{pmatrix}
w & Z & iz\\ X & A & -\Bbb IZ^*\\ ix & -X^*\Bbb I & -\bar w
\end{pmatrix}
$$
with blocks of sizes $1$, $n$, and $1$, $w\in\Bbb C$, $x,z\in\Bbb
R$, $X\in\Bbb C^n$, $Z\in\Bbb C^{n*}$, and $A\in\frak u(p,q)$ such
that $w-\bar w+\tr(A)=0$. Here $\Bbb I$ is the diagonal matrix with
the first $p$ entries equal to $1$ and the remaining $q$ entries equal
to $-1$. 

It is easy to show that the subgroup $Q\subset G$ corresponds to
matrices of the form
$$
\begin{pmatrix}
\ph & 0 & ia\ph \\ 0 & \Ph & 0\\ 0 & 0 & \bar\ph^{-1}  
\end{pmatrix},
$$
with $\ph\in\Bbb C\setminus\{0\}$, $a\in\Bbb R$ and $\Ph\in U(p,q)$
such that $\tfrac{\ph^2}{|\ph|^2}\det(\Ph)=1$. 

\subsection{}\label{5.2}
Next we need an analog of the pair $(i,\al)$ introduced in \ref{3.5}.
As before we start with a manifold $M$ of dimension $2n+1$, so again
$\tilde G=PGL(2n+2,\Bbb R)$. We will use a block decomposition into
blocks of sizes $1$, $1$, $n$, and $n$ as before.
The right choice turns out to be
\begin{gather*}
  i\begin{pmatrix}\ph & 0 & ia\ph \\ 0 & \Ph & 0\\ 0 & 0 & \bar\ph^{-1} 
  \end{pmatrix}:=
\begin{pmatrix}
|\ph| & -a|\ph| & 0 & 0\\ 0 & |\ph|^{-1} & 0 & 0\\ 0 & 0 &
 \Re(\tfrac{|\ph|}{\ph}\Ph) & -\Im(\tfrac{|\ph|}{\ph}\Ph) \\0 & 0 &
 \Im(\tfrac{|\ph|}{\ph}\Ph) & \Re(\tfrac{|\ph|}{\ph}\Ph)\end{pmatrix},\\
\al\begin{pmatrix} w & Z & iz\\ X & A & -\Bbb IZ^*\\ ix &
 -X^*\Bbb I & -\bar w\end{pmatrix}:=
\begin{pmatrix}
\Re(w) & -z & \Re(Z) & -\Im(Z)\\ x & -\Re(w) & -\Im(X^*\Bbb I) & -\Re(
X^*\Bbb I)\\ \Re(X) & \Im(\Bbb IZ^*) & \Re(A) & -\Im(A)+\Im(w)\\
\Im(X) & -\Re(\Bbb IZ^*) &  \Im(A)-\Im(w) & \Re(A)\end{pmatrix},
\end{gather*}
where $\Re$ and $\Im$ denote real and imaginary part, respectively,
and we write $\Im(w)$ for the appropriate multiple of the identity matrix.

There is an analog of Lemma \ref{3.6} (with similar proof), the only
change one has to make is that the map whose alternation has to be
used is given by 
$$
(X,Y,Z)\mapsto (\langle X_1,\Bbb IY_1\rangle+\langle X_2,\Bbb IY_2\rangle)
\begin{pmatrix}-Z_2 \\ Z_1\end{pmatrix}. 
$$
This map has similar properties as the one from \ref{3.6} so the
analogs of Theorem \ref{3.6} and Corollary \ref{3.6} hold. 

Concerning the structure of $\ker(\square)$ the situation is also
similar to the case of Lagrangean contact structures, since the
decomposition of $\ker(\square)$ can be determined from the
complexifications of $\fg$ and $\frak p$ which are the same in both
cases. The only difference is that the two irreducible components for
$n=1$ respectively the two irreducible components contained in
homogeneity $1$ in the case $n>1$ in the Lagrangean case correspond to
only one component here. This component however has a complex
structure and it consists of maps $\fg_{-1}\wedge\fg_{-2}\to\fg_1$
which are complex linear in the first variable respectively maps
$\La^2\fg_{-1}\to\fg_{-1}$, which are conjugate linear in both
variables. For $n>1$ this component is a torsion which is up to
a nonzero multiple given by the Nijenhuis tensor. Vanishing of this
component is equivalent to torsion freeness and to integrability of
the almost CR structure, see \cite[4.16]{C-S}. 

\begin{thm*}
Let $(M,H,J)$ be a partially integrable almost CR structure and let
$(p:\Cal G\to M,\om)$ be the corresponding regular normal parabolic
geometry of type $(G,P)$. Then the parabolic geometry $(\Cal
G\x_Q\tilde P\to\Cal P_0(TM),\tilde\om_\al)$ constructed using the
extension functor associated to the pair $(i,\al)$ from \ref{5.1} is
regular and normal if and only if $\om$ is torsion free, i.e.~the
almost CR structure is integrable. 
\end{thm*}
\begin{proof}
Apart from some numerical factors which cause no problems, this is
completely parallel to the proof of Theorem \ref{3.8}.  
\end{proof}

Hence the direct relation between the regular normal parabolic
geometries associated to a partially integrable almost CR structure
respectively to the associated path geometry of chains works exactly
on the the subclass of CR structures. 

\subsection{Applications}\label{5.3}
The developments of section \ref{4} can be applied to the CR case with
only minimal changes. In analog of Lemma \ref{4.3}, one obtains
$S\in\otimes^3 H_x^*\otimes H_x$, which is the complete symmetrization
of
$$
(\xi,\eta,\zeta)\mapsto \Cal L(\xi,J(\eta))J(\zeta),
$$
where $J$ is the almost complex structure on $H$. 

\begin{thm*}
Let $(M,H,J)$ be a CR structure.

\noindent
(1) There is no linear connection on $TM$ which has the chains among
its geodesics.

\noindent
(2) The path geometry of chains is torsion free if and only if the CR
structure is locally flat.

\noindent
(3) The almost complex structure $J$ can be reconstructed up to sign
from the harmonic curvature of the associated path geometry of chains.
\end{thm*}
\begin{proof}
  The only change compared to section \ref{4} is that one has to
  extend $S$ to the complexified bundle $H\otimes\Bbb C$. As in the
  proof of Theorem \ref{4.4} one then reconstructs the subset
  $H_x^{1,0}\cup H_x^{0,1}\subset H_x\otimes\Bbb C$ for each $x\in M$,
  i.e.~the union of the holomorphic and the anti--holomorphic part.
  This union determines $J$ up to sign.
\end{proof}
This theorem now also implies that the signature of the CR structure,
which is encoded in $\Cal L(-,J(-))$, can be reconstructed from the
path geometry of chains. As a corollary, we obtain a completely
independent proof of the analog of Corollary \ref{4.4}, which is due to
\cite{Cheng} for CR structures:
\begin{cor*}
A contact diffeomorphism between two CR manifolds which maps chains to
chains is either a CR isomorphism or a CR anti--isomorphism.  
\end{cor*}

%%%%%%%%%%%%%%%%%%%%%%%%%%

\end{document}